%%%%%%%%%%%%%%%%%%%%%%%%%%%%%%%%%%%%%%%%%%%%%%%%    
%
%        THIS IS A  PLAIN TeX FILE
%
%%%%%%%%%%%%%%%%%%%%%%%%%%%%%%%%%%%%%%%%%%%%%%%%

\magnification=1200

\font\titfont=cmr10 at 12 pt

\font\headfont=cmr10 at 12 pt

%\font\AAA=Times.dfont  at 12pt
 %\font\BBB=Times.dfont at 8pt

%\font\AAA=cmr10 at 12pt
%\font\BBB=cmr10 at 8pt

\overfullrule=0in

\def\boxit#1{\hbox{\vrule
 \vtop{%
  \vbox{\hrule\kern 2pt %
     \hbox{\kern 2pt #1\kern 2pt}}%
   \kern 2pt \hrule }%
  \vrule}}

  \def\harr#1#2{\ \smash{\mathop{\hbox to .3in{\rightarrowfill}}\limits^{\scriptstyle#1}_{\scriptstyle#2}}\ }

\def\IFF{{\qquad \iff \qquad}}
\def\GL{{\rm GL}}
\def\Symr{{\rm Sym}^2_{\bbr}}
\def\HSym{{\rm HSym}^2}
\def\HSymr#1{{\rm HSym}^2_{\bbr}(#1)}
\def\bra#1#2{\langle#1, #2\rangle}
\def\ss{\subset}
\def\half{\hbox{${1\over 2}$}}
\def\smfrac#1#2{\hbox{${#1\over #2}$}}

\def\dist{{\rm dist}}

\def\Hess{{\rm Hess}}

\def\tr{{\rm tr}}

\def\Hom{{\rm Hom\,}}
\def\det{{\rm det}}
\def\End{{\rm End}}
\def\Sym{{\rm Sym}^2}

\def\arr{\longrightarrow}

\def\rn{\bbr^n}

\def\Int{{\rm Int}}

\def\Symn{{\Sym(\rn)}}

\def\Theorem#1{\medskip\noindent {\bf THEOREM \bf #1.}}
\def\Prop#1{\medskip\noindent {\bf Proposition #1.}}
\def\Cor#1{\medskip\noindent {\bf Corollary #1.}}
\def\Lemma#1{\medskip\noindent {\bf Lemma #1.}}
\def\Remark#1{\medskip\noindent {\bf Remark #1.}}

\def\Def#1{\medskip\noindent {\bf Definition #1.}}

\def\Ex#1{\medskip\noindent {\bf Example \bf    #1.}}

\def\pf{\medskip\noindent {\bf Proof.}\ }
\def\qed{\hfill  $\vrule width5pt height5pt depth0pt$}
\def\equdef{\buildrel {\rm def} \over  =}

\def\mathqed{\hfill  \vrule width5pt height5pt depth0pt}

\def\n{\nabla}

      \def\cb{{\cal B}}  
   \def\cp{{\cal P}}   
   
\def\ce{{\cal E}}   
\def\ch{{\cal H}}   
   
\def\cd{{\cal D}}

\def\cp{{\cal P}}
\def\cf{{\cal F}}

\def\vf{\varphi}

\def\wt{\widetilde}

\def\and{\qquad {\rm and} \qquad}
\def\arr{\longrightarrow}
\def\ol{\overline}
\def\bbr{{\bf R}}
\def\bbc{{\bf C}}

\def\bbj{{\bf J}}

\def\a{\alpha}
\def\b{\beta}

\def\e{\epsilon}

\def\l{\lambda}
\def\o{\omega}

\def\s{\sigma}

\def\z{\zeta}

\def\D{\Delta}
\def\L{\Lambda}
\def\G{\Gamma}
\def\O{\Omega}

\def\psh{plurisubharmonic }

\def\lloc{L^1_{\rm loc}}
\def\dbar{\ol{\partial}}

\def\bo{\partial \Omega}

\def\PSH{{ \rm PSH}}
\def\SH{{\rm SH}}

\def\Symn{\Sym(\rn)}
 
\def\USC{{\rm USC}}
\def\fa{{\rm\ \  for\ all\ }}

\def\ob{\overline{\O}}

\def\AA{1}
\def\BB{2}
\def\CC{3}

\def\EE{4}
\def\FF{5}
\def\II{6}
\def\HH{7}
\def\GG{9}

\def\JJ{8}

 \def\bbf{{\bf F}}

\def\HB{{\rm H\Sym}}
\def\FPSH{$F(J)$-plurisubharmonic }   \def\FJ{F(J)} \def\FxJ{F_x(J)}
\def\FPSHP#1{$F(#1)$-plurisubharmonic } 
\def\FJP#1{F(#1)}
\def\FJXP#1#2{F_{#1}(#2)} 

\centerline{\titfont POTENTIAL THEORY ON ALMOST COMPLEX MANIFOLDS }
\medskip
\bigskip

\centerline{\titfont F. Reese Harvey and H. Blaine Lawson, Jr.$^*$}
\vglue .9cm
\smallbreak\footnote{}{ $ {} \sp{ *}{\rm Partially}$  supported by
the N.S.F. }

 \vskip .3in
\centerline{\bf ABSTRACT} \bigskip
  \font\abstractfont=cmr10 at 10 pt
{{\parindent= .53in
\narrower\abstractfont \noindent
Pseudo-holomorphic curves  on almost complex manifolds have been much more intensely  studied
than their  ``dual'' objects, the  plurisubharmonic functions.  These
functions are  standardly   defined by requiring that the restriction to each pseudo-holomorphic curve
be subharmonic.  In this paper subharmonic functions are defined
by applying  the viscosity approach   to a version of the complex hessian which
exists intrinsically  on any almost  complex manifold. 
Three theorems are proven. The first is a restriction theorem which establishes the 
equivalence of our definition with the ``standard''  definition.
In the second theorem, using our ``viscosity'' definitions, 
 the Dirichlet problem is solved for the 
 complex Monge-Amp\`ere equation in both the homogeneous and
 inhomogeneous forms.
These  two  results are based on theorems found in [HL$_3$] and [HL$_2$] respectively. 
Finally, it is shown that the plurisubharmonic functions considered here agree with
the plurisubharmonic distributions. In particular, this proves a conjecture of Nefton Pali.
\vskip .1in

}}

\vskip .3in

\centerline{\bf TABLE OF CONTENTS} \bigskip

{{\parindent= .1in\narrower\abstractfont \noindent

\qquad 1. Preliminaries: Almost Complex Structures and 2-Jets. \smallskip

\qquad 2.      The Complex Hessian.
\smallskip

\qquad 3.  \FPSH Functions. \smallskip

\smallskip

\qquad 4.  A Local Coordinate Expression for the Real Form of the Complex Hessian.

\smallskip

\qquad 5.   Restriction of \FPSH Functions.
\smallskip

\qquad 6.   The Equivalence of \FPSH  Functions 
and 

\qquad\ \ \  \  Standard Plurisubharmonic Functions.
\smallskip

\qquad 7.   The Dirichlet Problem.
\smallskip

\qquad 8.  Distributionally Plurisubharmonic Functions.
\smallskip

\qquad 9.   The Non-Equivalence of Hermitian and Standard Purisubharmonic Functions.
\smallskip

\qquad\qquad Appendix A. The Equivalence of the Various Notions of Subharmonicity

\qquad\qquad \qquad\qquad\ \ \   for Reduced Linear Elliptic Equations.

}}

\vfill\eject

 %%%%%%%%%%%%%%%%%%%%%%%%%%%%%%%%%%%%%%%%%%%%%%%%%
  %%%%%%%%%%%%%%%%%%%%%%%%%%%%%%%%%%%%%%%%%%%%%%%%%
  %%%%%%%%%%%%%%%%%%%%%%%%%%%%%%%%%%%%%%%%%%%%%%%%%
  %%%%%%%%%%%%%%%%%%%%%%%%%%%%%%%%%%%%%%%%%%%%%%%%%

\noindent{\headfont 0.\  Introduction.}
\medskip
The purpose of this paper is to develop an intrinsic potential theory on a general almost
complex manifold $(X,J)$.  Our methods are based  on results established in [HL$_{2,3,4}$].
% and [HL$_2$].
In particular, we  study an extension, to this general situation, of the classical notion
of a plurisubharmonic function.  For smooth  functions  $\vf$ many equivalent definitions
of plurisubharmonicity are available.  Several are given in Section 2.  For instance, one can
require that $\, i\partial\dbar \vf \geq0$, or $\ch(\vf)\geq0$ where the complex hessian
$\ch(\vf)(V,W) \equiv (\half V\overline W +\half \overline W  V -{i\over 2}J[V,\overline W])(\vf)$
is a bilinear form on $T_{1,0}X$.  Its real form $H\equiv {\rm Re}\,\ch$
is computed in Section 4 to be $H(\vf)(v,v) = (vv+(Jv)(Jv) + J[v,Jv])(\vf)$, yielding a third 
equivalent condition $H(\vf)\geq0$.

At any given point $x\in X$, each of these equivalent definitions only 
depends on the reduced 2-jet of $\vf$ at $x$.  Consequently, one obtains a subset 
$\FJ \ss J^2(X)$ of the 2-jet bundle of $X$, which consists
 of the $J$-plurisubharmonic jets.  
 
 For an upper semi-continuous function $u$
on $X$, the previous definitions cannot be applied directly to $u$.
However, they can be applied to a smooth ``test function'' $\vf$ for $u$ at $x$.
(See Section 3 for more details.) This yields our general definition of $\FJ$-plurisubharmonic
functions.

Our first main result is the Restriction Theorem \FF.2 which states that: {\bf  if $u$
is  $\FJ$-plurisubharmonic and $(X',J')$ is an almost complex submanifold of $(X,J)$,
then the restriction of $u$ to $X'$ is  $F(J')$-plurisubharmonic.}

An upper semi-continuous function $u$ on $(X,J)$, with the property that its restriction 
to each pseudo-holomorphic curve in $X$ is subharmonic, will be called 
{\sl standardly  plurisubharmonic} (Definition \II.1).
On such curves the complex structure is integrable and all of the many
definitions of subharmonic are well known to agree. Consequently, the special case of our 
restriction theorem, where $X'$ has dimension one, states that each $\FJ$-plurisubharmonic
function is  plurisubharmonic in the standard sense.  
The converse is also true due to the classical theorem
of Nijenhuis and Woolf, that there exist
pseudo-holomorphic curves through any point of $X$ with any prescribed tangent line.
This gives our first main result, the equivalence  of the two notions of plurisubharmonicity (Theorem \II.2).

%There is another notion of plurisubharmonicity for upper semi-continuous functions
%which uses the family of $F(J)$-{\sl harmonic functions} ( in  [HL$_2$]-terminology --
%also called  {\sl maximal} plurisubharmonic functions in [BT]-terminology 
%or  $F(J)$-{\sl solutions} in the viscosity world).  By definition $u$ is plurisubharmonic in
%this sense if $u-h$ satisfies the maximum principle for all $F(J)$-harmonics $h$ (See Definition \HH.%10).
%The equivalence of this notion
%with those above is established in Theorem \HH.11.

The second main result,  Theorem \HH.5, solves the {\bf Dirichlet problem
for  ${\bf F(J)}$-harmonic functions} in the general almost complex setting.
Here the usefulness of the $\FJ$-notion over the standard notion is striking.
The result  applies simultaneously to both the homogeneous and the 
 inhomogeneous equations, and the proof rests on  basic results in  [HL$_2$]. 
The main statement is essentially the following. Fix a real volume form $\l$ on $X$ 
which is orientation compatible with  $J$.  
Let $\O\ss X$ be a compact domain with smooth boundary $\bo$.  Fix a  continuous
functions $f\in C(\ob)$ with $f\geq0$. Given $\vf\in C(\bo)$, the  Dirichlet Problem asks
for existence and uniqueness of $J$-plurisubharmonic viscosity solutions $u\in C(\ob)$
to the equation 
$$
(i\partial\dbar u)^n \ =\ f \l.
$$
with boundary values $u\bigr|_{\bo} =\vf$. Existence and Uniqueness are established whenever $\O$ has
a strictly $J$-plurisubharmonic defining function. See Theorem \HH.4 for full details.

At the end of Section \HH \ we discuss some other   equivalent notions of 
$F(J)$-harmonic functions.

The hard work in proving the restriction theorem occurs in [HL$_3$] while the
hard work in solving the Dirichlet problem occurs in [HL$_2$]. 
 In both cases, the key to being able to apply these results to almost complex 
 manifolds is the local coordinate expression given in Proposition \EE.5,
 which states that the subequation $\FJ$ is locally jet equivalent to the standard constant
 coefficient case on local charts in $\bbr^{2n}=\bbc^n$ equipped with a standard complex structure.
 The methods must be adapted somewhat unless $f\equiv0$ or $f>0$ (see Section \HH).

On any almost  complex manifold there is also the notion of a plurisubharmonic distribution.
It has been proved by Nefton Pali [P] that  if $u$ is  \psh in the standard sense,
then $u$ is in $\lloc(X)$ and 
defines a \psh distribution.  Pali conjectures that the converse is also true and proves a  
partial result in this direction.  In Section 8 we prove the full conjecture, thereby
showing that  on any almost complex
manifold all three notions of plurisubharmonicity (standard, viscosity and distributional) are equivalent.

The argument given in Section 8 reduces this nonlinear result to a corresponding result
for linear elliptic subequations.  This technique applies  to any convex subequation which
 is ``second-order complete'' in a certain precise sense.
In an appendix, linear elliptic operators  $L$  with smooth coefficients and no zero-order term
are discussed.  The notions of 
a viscosity $L$-subharmonic function and an $L$-subharmonic distribution are shown to
be equivalent.   First we show that an upper semi-continuous function  is   $L$-subharmonic in the viscosity sense if and only if it is  $L$-subharmonic in the ``classical''  sense 
 of being ``sub-the-$L$-harmonics''.  Such functions
are known to be $\lloc$ and to give an $L$-subharmonic distribution (Theorem A.6).  Conversely, any
$L$-subharmonic distribution has a unique upper semi-continuous, $\lloc$ representative
which is ``classically'' subharmonic.  It is  important not only  that this representative be unique but also that it can be 
obtained canonically from the $\lloc$-class of the function by taking the 
essential upper semi-continuous regularization $\wt u(x) = {\rm ess} \, \limsup_{y\to x} u(y)$.
Thus the choice of upper semi-continuous function in the $\lloc$-class is independent of 
the operator $L$.  (It is actually  the smallest upper semi-continuous  representative.)
This fact is essential for the arguments of Section 8. 

The proof  for this linear case combines techniques from  distributional potential theory, 
classical potential theory and viscosity  potential theory. 
It appears that there is no specific reference for this particular result, so
we have included an appendix which outlines the theory.
(Also see Section 9 in [HL$_4$].)

\medskip
 
\noindent
{\bf Note 1.} 
We note that when $(X,J)$ is given a hermitian metric,   there is a notion of   ``metric'' plurisubharmonic
functions, defined via the hermitian symmetric part of the riemannian hessian
(cf. [HL$_2$]).  For metrics with the property that every holomorphic curve in $(X,J)$ is  minimal
 (e.g., if the  K\"aher form  is closed), these ``metric'' plurisubharmonic
functions are standard (Theorem 9.3).  However in Section 9 examples are given which,
in the general case,  strongly differentiate the metric notion
from the intrinsic    notion of plurisubharmonicity studied in this paper.

\medskip

\noindent
{\bf Note 2.}  In the first appearance of this article we solved the Dirichlet problem 
for  $f\equiv0$ (and the methods also  applied to $f>0$).  Shortly afterward
  Szymon Pli\'s posted a paper which also  studied the Dirichlet problem on almost complex
manifolds [Pl].  His result, which is the analogue in the almost complex case of a main 
result in [CKNS],  is quite different from ours.  He considered the case  $f>0$,  
 assumed all data   to be smooth,  and established complete regularity of the solution. Our
 result holds for arbitrary continuous boundary data, and interior regularity is known to fail.
 
 Pl\'is subsequently posted a second  version of [Pl] which offered an alternative proof
 of our result by taking a certain limit of his solutions.
 He  also established Lipschitz continuity of the solution under restrictive 
 hypotheses on the boundary data. In the current version of this paper we allow $f\geq0$.

\medskip

It is assumed throughout that $J$ is  of class $C^\infty$, however  for all results
which do not involve distributions, $C^2$ is sufficient.

The authors are indebted to Nefton Pali and Jean-Pierre Demailly
for useful conversations and comments regarding this paper.
We would also  like to thank Mike Crandall, Hitoshi Ishii, Andrzej Swiech and Craig Evans 
for helpful discussions regarding the Appendix.

\vfill\eject

 %%%%%%%%%%%%%%%%%%%%%%%%%%%%%%%%%%%%%%%%%%%%%%%%%
  %%%%%%%%%%%%%%%%%%%%%%%%%%%%%%%%%%%%%%%%%%%%%%%%%
  %%%%%%%%%%%%%%%%%%%%%%%%%%%%%%%%%%%%%%%%%%%%%%%%%
  %%%%%%%%%%%%%%%%%%%%%%%%%%%%%%%%%%%%%%%%%%%%%%%%%

\noindent{\headfont \AA.\  Preliminaries}
\medskip
\centerline{\bf Almost  Complex Structures}
\medskip
Let $(X,J)$ be an almost complex manifold.  Then there is a  natural decomposition
$$
TX\otimes_\bbr \bbc \ =\ T_{1,0} \oplus T_{0,1} 
$$
where    $T_{1,0}$ and $T_{0,1}$ are the  $i$ and $-i$-eigenspaces of $J$ respectively. 
The projections of $TX\otimes_\bbr \bbc$ onto these complex subspaces are given explicitly
by $\pi_{1,0} = \half(I-iJ)$ and $\pi_{0,1} = \half(I+iJ)$.
 
 The dual action of $J$ on $T^*X$ shall be denoted by $J$ as well.
(For $\a\in T^*_xX$, $(J\a)(V) \equiv  \a(JV)$ for $V\in T_xX$.)
Then again we have a  natural decomposition
$$
T^*X\otimes_\bbr \bbc \ =\ T^{1,0} \oplus T^{0,1} 
$$
into the  $i$ and $-i$-eigenspaces of $J$ respectively with 
  projections again  given explicitly
by $\pi^{1,0} = \half(I-iJ)$ and $\pi^{0,1} = \half(I+iJ)$.
Moreover, we have
$$
(T_{1,0})^* \ =\ T^{1,0} 
\and
(T_{0,1} )^*\ =\ T^{0,1} 
$$

There is also a bundle splitting
$$
\L^kT^*X \otimes_\bbr \bbc  \ =\ \bigoplus_{p+q=k} \L^{p,q}  X 
\eqno{(\AA.1)}
$$
where $\L^{p,q} X $ is the $i(p-q)$-eigenspace of $J$ acting as a derivation on 
$\L^*T^*X \otimes_\bbr \bbc$.  Let $\ce^{p,q} X $ denote the smooth sections of $\L^{p,q} X $.
Then there are natural operators
$$
\partial : \ce^{p,q}(X) \ \to\  \ce^{p+1,q}(X)
\and
\dbar : \ce^{p,q}(X) \ \to\  \ce^{p,q+1}(X)
$$
defined by restriction and projection of the exterior derivative $d$.

Under complex conjugation one has that
$$
\overline{{\L^{p,q}} X } \ =\ \L^{q,p} X 
$$
and in particular each $\L^{p,p} X $ is conjugation-invariant and 
decomposes into a real and imaginary part.  Let $\L^{p,p}_\bbr X $ denote the real 
part.

A  {\sl pseudo-holomorphic map} $\Phi:(X',J')\to (X,J)$ between almost 
complex manifolds is a smooth
map whose differential $\Phi_*$  satisfies  $\Phi_* J'  = J  \Phi_*$ at every point.
Thus the pull-back $\Phi^*  :\L^*T^*X\to  \L^*T^*X'$ is also 
compatible with the almost complex structures (acting as derivations)
and therefore preserves the bigrading in  (\AA.1).
It follows that:
$$
 {\rm  The\  operators  \ }  \partial\ {\rm and}\ 
  \dbar \ {\rm commute\ with\ the\  pullback\ }  \Phi^* \ {\rm  on\ smooth\ forms.}
\eqno{(\AA.2)}
$$

%%%%%%%%%%%%%%%%%%%%%%%%%%%%%%%%%%%%%%%%%%%%%%%%
\vfill\eject
\centerline{\bf 2-Jets}
\medskip

Denote by $\overline J^2(X)$ the vector bundle of 2-jets on an arbitrary smooth manifold $X$,
and let $\overline J_x(u)\in \overline J_x^2(X)$ denote the 2-jet of a smooth function $u$ at $x$.  
The bundle ${J}^2(X)$ of reduced 2-jets is defined to be the quotient
of $\overline{J}^2(X)$ by the trivial line bundle corresponding to the value of the
function $u$. We will be primarily interested in the space of reduced 2-jets in this paper
(and so we have chosen a simpler notation for it).
The bundle $\Sym(T^*X)$ of quadratic forms on $TX$ has a natural embedding as a subbundle
of  ${J}^2(X)$ via the hessian at critical points.
  Namely, if $u$ is a function with a critical point at $x$, then
$V(W(u)) = W(V(u)) + [V,W](u) = W(V(u))$ is a well defined symmetric bilinear form
$(\Hess_x u) (V,W)$ on $T_xX$ for arbitrary vector fields $V,W$ defined near $x$.

Thus  there is a short exact sequence of bundles
$$
0\ \ \arr\ \ \Sym(T^*X)\ \ \arr\ \  {J}^2(X) \ \ \arr\ \ T^*X\ \ \arr\ \ 0.
\eqno{(\AA.3)}
$$
However, for functions $u$ with $(du)_x\neq0$, there is no natural definition
of $\Hess_x u$, i.e., the sequence (\AA.3) has no natural splitting.
There is a natural cone bundle $\cp(X)\ss \Sym(T^*X) \ss  {J}^2(X)$
defined by
$$
\cp_x(X)\ \equiv \ \{H\in   \Sym(T_x^*X): H\geq0\}\ \cong\ \{ {J}_x(u) : u(x) =0 \ {\rm and\ } u\geq0 \ {\rm near\ } x\}
\eqno{(\AA.4)}
$$

In local coordinates $x\in \bbr^N$ for $X$, the first and second derivatives comprise the reduced 
2-jet ${J}_x(u)$ of a function $u$.  That is,
$$
{J}_x(u) \ \cong\ (D_xu, D^2_xu)\ \in\ \bbr^N \times \Sym(\bbr^N) \ \equdef 
\  \bbj^2(\bbr^N) \ \equiv \   \bbj^2
\eqno{(\AA.5)}
$$
where $D_x u \equiv ({\partial u\over \partial x_1}(x),...,{\partial u\over \partial x_N}(x))$ and
$D^2_xu \equiv \left(   {\partial^2 u\over \partial x_i\partial x_j} (x)   \right)$.

The isomorphism (\AA.5) says that ${J}_x^2(X)  \cong \bbr^N\times\Sym(\bbr^N)
=  \bbj^2(\bbr^N)$.
The standard notation $(p,A) \in \bbr^N\times \Sym(\bbr^N)$
 will be used for coordinates $(p,A)= (Du, D^2u)$ on $ \bbj^2(\bbr^N)$.

Using these coordinates, we have
$$
\Sym_x(T^*X)\ \cong\ \{(0,A) : A\in \Sym(\bbr^N)\}\ \ss\  \bbj^2
\eqno{(\AA.6)}
$$
and
$$
\cp_x(X) \ \cong\ \{(0,A) : A \geq 0\}\ \ss\  \bbj^2
\eqno{(\AA.7)}
$$

%%%%%%%%%%%%%%%%%%%%%%%%%%%%%%%%%%%%%%%%%%%%%%%
%%%%%%%%%%%%%%%%%%%%%%%%%%%%%%%%%%%%%%%%%%%%%%%
%%%%%%%%%%%%%%%%%%%%%%%%%%%%%%%%%%%%%%%%%%%%%%%
%%%%%%%%%%%%%%%%%%%%%%%%%%%%%%%%%%%%%%%%%%%%%%%
%%%%%%%%%%%%%%%%%%%%%%%%%%%%%%%%%%%%%%%%%%%%%%%

\vfill \eject
%\vskip .3in
\noindent{\headfont  \BB. The Complex Hessian}\medskip

Suppose $(X,J)$ is an almost complex manifold.
Note that for any real-valued $C^2$-function $u$ we have 
$
d \dbar u = \dbar^2 u + \partial \dbar u + \e^{2,0}
$
where $\e^{2,0}$ is a continuous (2,0)-form. Taking the complex conjugate
and adding terms gives $d d u = 0$,  which, after taking $(2,0)$- and then $(1,1)$-components,
gives us that $\e^{2,0} = -\partial^2 u$ and 
$$
\partial \dbar u \ =\ -\dbar \partial u.
$$

\Def{\BB.1} The {\bf $i \partial\dbar$-hessian} of a smooth 
real-valued function $u$ on $X$ is obtained by applying 
 the intrinsically defined real operator
$$
 i \partial \dbar : \ C^\infty(X) \ \arr\  \G(X, \L_\bbr^{1,1} X ),
\eqno{(\BB.1)}
$$

From (\AA.2) we conclude the following.

\Prop{\BB.2} {\sl The $i \partial\dbar$-hessian operator commutes with pull-back under any 
pseudo-holomorphic map between almost complex manifolds.
}\medskip

There is a second useful way of looking at $i\partial \dbar u$.
Recall that a (complex-valued) bilinear form $\cb$ on a complex vector space is {\bf hermitian} if
$\cb(V,W)$ is complex linear in the first variable and  complex anti-linear in the second variable. It is
further called {\bf symmetric} if it
satisfies
$$
\cb(V,W)\ =\ \overline{ \cb(W,V)}.
$$
When, in addition, $\cb(V,V)\geq0$ for all $V$, 
we say $\cb$ is {\bf  positive} and denote this by $\cb\geq0$.
Let $\HB(T_{1,0}X)$ denote the bundle of hermitian symmetric forms on $T_{1,0}X$. 

Recall the bundle  isomorphism 
$$
\L^{1,1}_\bbr  X  \ \cong\ \HB(T_{1,0}X)
\eqno{(\BB.2)}
$$
 sending 
$$
\b \in \L_{\bbr}^{1,1}X\qquad {\rm to}\qquad \cb(V,W) \equiv -i\b(V,\overline W).
$$
We say $\b\geq0$ if $\cb\geq0$.

\Def{\BB.3}  Under the isomorphism (\BB.2) the $i \partial\dbar$-hessian of $u$ 
becomes the {\bf complex hessian $\ch(u)$ of $u$}.  Namely 
$$
\ch(u)(V,W) \ =\  (\partial\dbar u ) (V,\overline W)
\eqno{(\BB.3)}
$$
is a section of $\HSym(T_{1,0}X)$, and 
$$
\ch : C^\infty(X) \ \arr\ \G(X, \HSym(T_{1,0}X))
$$
will be called the {\bf complex hessian operator}.

\medskip

The complex hessian can be computed in various ways. The main formula we employ is the following.

\Prop{\BB.4}  {\sl For $u\in C^\infty(X)$ and  each $V,W \in \G(X, T_{1,0})$,}
$$
\ch(u)(V,W) \ =\  \left( \half V\overline W + \half \overline W V - \smfrac i 2 J[V, \overline W] \right)(u)
\eqno{(\BB.4)}
$$

\Remark{\BB.5} We can express (\BB.4) more succinctly by saying that:  
for each $V,W \in \G(X, T_{1,0})$, $\ch(V,W)$ is the second-order scalar differential operator
$$
\ch (V,W) \ =\  \half V\overline W + \half \overline W V - \smfrac i 2 J[V, \overline W] 
\eqno{(\BB.4)'}
$$

\pf
First recall that for arbitrary sections  $V$ and $W$   of $TX \otimes_\bbr \bbc$ and any complex 1-form $\a$,
 the exterior derivative of $\a$  satisfies:
$$
(d\a) \left( V,\overline W  \right)\ = \  V \left(  \a(\overline W)  \right)    -
  \overline W \left(  \a(V)  \right)  -\a \left( [V, \overline W]  \right).
$$
Now assume that $V$ and $W$ are both of type $1,0$ and 
take $\a = \dbar u$.  Then  $\a(\overline W) = \overline W (u)$
while $\a(V)=0$.  Since $(\partial \dbar u) (V, \overline W) = (d\a)(V, \overline W)$, we have
$$
(\partial\dbar u)(V, \overline W)\ =\ 
V  \left( \overline W ( u) \right) - \dbar u \left( [V, \overline W]\right ) \ =\ 
V  \left( \overline W ( u) \right)  -  [V, \overline W]^{0,1} (u).
\eqno{(\BB.5)}
$$
Take $\a = - \partial u$ and note that $\a(\overline W) = 0$
while $\a(V)=  - V(u)$.  
 Since $ - (\partial \dbar u) (V, \overline W) =  - (d\a)(V, \overline W)$, we have
$$
 - (\dbar \partial u)(V, \overline W)\ =\ 
 \overline W  \left ( V (u) \right ) + \partial u \left( [V, \overline W]\right) \ =\ 
 \overline W  \left ( V (u) \right )  + [V, \overline W]^{1,0} (u).
\eqno{(\BB.6)}
$$
Finally, using $J du  =  i(\partial u-\dbar u)$, 
or equivalently that $J=i$ on $T_{1,0}X$ and $J=-i$ on $T_{0,1}X$,
we see that the average of these last  two formulas for $\ch(u)(V, \overline W)$ is 
given by (\BB.4). \qed

%%%%%%%%%%%%%%%%%%%%%%%%%%%%%%%%%%%%%%%%%%%%%%%
%%%%%%%%%%%%%%%%%%%%%%%%%%%%%%%%%%%%%%%%%%%%%%%
%%%%%%%%%%%%%%%%%%%%%%%%%%%%%%%%%%%%%%%%%%%%%%%
%%%%%%%%%%%%%%%%%%%%%%%%%%%%%%%%%%%%%%%%%%%%%%%
%%%%%%%%%%%%%%%%%%%%%%%%%%%%%%%%%%%%%%%%%%%%%%%

\vskip .3in
%\vfill\eject
\noindent{\headfont  \CC. \FPSH Functions.}\medskip

It is natural, and useful, to have a definition of plurisubhamonic functions on $(X,J)$ expressed
purely in terms of the 2-jets of those functions. Specifically, one would like such functions, when smooth,
to be defined by constraining their 2-jets to a subset $\FJ \ss J^2(X)$ of the 2-jet bundle, 
and then pass to general  upper semi-continuous functions by viscosity techniques.
In this section we  give such a definition  using the complex hessian $\ch$.  
First, for smooth functions the concept is straightforward.

\Def{\CC.1}  A smooth real-valued function $u$ on $X$ is called  {\bf \FPSH} 
if $\ch(u)\geq0$, i.e., the hermitian symmetric bilinear form $\ch_x(u)$ 
is positive semi-definite at all points  $x\in X$.  
Moreover, if $\ch_x(u) > 0$ is positive definite at  all points $x\in X$, we say that $u$ is {\bf strictly} 
\FPSH.
\medskip

Proposition \BB.4 implies that   at a point $x\in X$, $\ch_x(u)$
   depends only on ${J}_x(u)$, 
the reduced 2-jet of    $u$ at $x$.   In particular, the condition $\ch_x(u)\geq 0$
(equivalently $i\partial\dbar u\geq0$) depends 
only on the jet ${J}_x(u)$ of $u$ at $x$.
 Hence we can define  plurisubharmonicity 
for a jet $\bbj \in  {J}^2_x(X)$ as follows.

\Def{\CC.2}  A jet $\bbj \in   {J}^2_x(X)$  is said to be   {\bf \FPSH} if for any smooth
function $u$ with ${J}_x(u)= \bbj$, we have 
$$
\ch_x(u) \geq0.
\eqno{(\CC.1)}
$$
The set of \FPSH jets on $X$  will be denoted by $\FJ$.

\medskip

We are now in a position to  broaden the notion of $\FJ$-plurisubharmonicity to the level 
of generality encountered in 
classical complex function theory.   Namely we consider functions 
 $u\in  \USC(X)$, the space of upper semi-continuous
functions on $X$ taking values in $[-\infty,\infty)$.  Take $u\in \USC(X)$ and fix $x\in X$.  A function $\vf$ which is $C^2$ in a neighborhood of $x$ 
is called a {\bf test function for $u$ at $x$} if $u-\vf\leq0$ near $x$ and $u-\vf=0$ at $x$.

\Def{\CC.3}  A function $u\in \USC(X)$ is {\bf \FPSH} if for each $x\in X$ and each
test function $\vf$ for $u$ at $x$, one has
$$
\ch_x(\vf) \geq0, \ \ {\rm i.e.\ \ }J_x\vf\ \in\  \FxJ.
$$
 
 Note that $u\equiv-\infty$ is \FPSH since there exist no test functions for $u$ at any point.

Definition \CC.3 should be an extension of Definition \CC.1 when $u$ is smooth.
For this to be true the following {\bf Positivity Condition} for $F$ (where $\cp(X)$ is defined by  (\AA.4))
$$
F+\cp(X)\ \ss\ F
\eqno{(P)}
$$
 must be satisfied for $F=\FJ$. (See Proposition 2.3 in [HL$_2$] for the details.)
 
 \Prop{\CC.4}  {\sl
 Each fibre $\FxJ$ of $\FJ$ is a convex cone, with vertex at the origin, containing
 the convex cone $\cp_x(X)$.  In particular, $\FJ$ satisfies  the Positivity Condition (P).
 }

\pf  It is easy to see that each fibre $\FxJ$ is a convex cone with vertex at the origin in
${J}_x^2(X)$ since $\ch_x(u)$ is linear in $u$.  It remains to show that $\FxJ$ contains
$\cp_x(X)$ as defined by (\AA.4).
Recall that each vector field $V$ of type $1,0$ is of the form $V=v-iJv$ where $v$ is a real vector field.
If $x$ is a critical point of $\vf$, then it is easy to compute from (\BB.4) that at $x$
$$
\ch (\vf) (V,V) \ =\ v(v(\vf)) + (Jv)(Jv(\vf))
\eqno{(\CC.2)}
$$
for all such  $1,0$ vector fields $V=v-iJv$. 
Now suppose that $\vf(x)=0$ and that $\vf\geq0$ near the point $x$.  Then by elementary calculus
$$
(d\vf)(x)\ =\ 0\and v(v(\vf))(x)\ \geq\ 0 \quad\forall\, v\in\G(X,TX).
\eqno{(\CC.3)}
$$
Combining (\CC.2) and (\CC.3)  proves that $\cp_x(X)\ss F_x(J)$.\qed

\Def{\CC.5}  
 A subset $F\ss   J^2(X) $  which satisfies both  the Positivity Condition
(P) and the {\bf  Topological Condition:}
$$
(i)\ \ F\ =\ \overline{\Int F}\qquad (ii) \ \ F_x\ =\ \overline{\Int F_x}\qquad 
(iii)\ \ \Int F_x \ =\  (\Int F)_x
\eqno{(T)}
$$
is called a  {\bf subequation} (cf.  [HL$_2$, Def. 3.9]).  
This condition (T) for $\FJ$ is a   consequence
of a jet equivalence for the complex hessian which is given in the next section.

\Def{\CC.6}  A function $u\in\USC(X)$ is called $F$-{\bf subharmonic} on $X$ if 
for each $x\in X$ and each test function $\vf$ for $u$ at $x$, we have $J^2_x(\vf)\in F_x$.
The set of such functions is denoted $F(X)$.\medskip

Starting with classical potential theory on $\bbc^n$ as motivation, 
many of the important  results concerning  plurisubhamonic functions on $\bbc^n$ 
 were extended  to general constant coefficient
subequations in euclidean space in one of our first papers on the subject [HL$_1$]. 
Most of these results were then generalized to any subequation $F$ on a manifold  [HL$_2$].
One can view this paper as coming full circle back to the complex setting by using viscosity methods to prove a new result in the almost complex case.

%%%%%%%%%%%%%%%%%%%%%%%%%%%%%%%%%%%%%%%%%%%%%%%
%%%%%%%%%%%%%%%%%%%%%%%%%%%%%%%%%%%%%%%%%%%%%%%
%%%%%%%%%%%%%%%%%%%%%%%%%%%%%%%%%%%%%%%%%%%%%%%
%%%%%%%%%%%%%%%%%%%%%%%%%%%%%%%%%%%%%%%%%%%%%%%
%%%%%%%%%%%%%%%%%%%%%%%%%%%%%%%%%%%%%%%%%%%%%%%
\vfill\eject
%\vskip .3in

	\noindent{\headfont  \EE.  A  Local Coordinate Expression for the Real Form of the Complex Hessian.}\medskip

The point of this section is to establish a formula for the complex hessian in a real coordinate
system on $X$. We begin by reviewing some standard algebra.
The space $\HSym(T_{1,0})$ of hermitian symmetric bilinear forms on $T_{1,0}$ has an alternate
description.  Recall  the standard isomorphism on complex vector spaces
$$
(T,J) \ \cong\ (T_{1,0}, i)
\eqno{(\EE.1)}
$$
given by mapping a real tangent vector $v\in T$ to $V=\half(v-iJv)$ with inverse 
$v=2{\rm Re}\,V$.

A real symmetric bilinear form $B\in \Sym_\bbr(T)$ is said to be {\bf hermitian}
(or $J$-{\bf hermitian})  if $B(Jv,Jv)=B(v,v)$ for all $v\in T$.
Let $\HSymr T$ denote the subspace of $\Sym_\bbr(T)$ consisting of $J$-hermitian forms.
Now (\EE.1) induces a (renormalized) isomorphism
$$
\HSymr T \ \cong\ \HSym(T_{1,0}),
\eqno{(\EE.2)}
$$
given by mapping ${\cal B} \in  \HSym(T_{1,0})$ to its {\bf real form}
$B\in \HSymr T $ defined by 
$$
B(v,w) \ \equiv \ {\rm Re} \, {\cal B}(v-iJv, w-iJw).
\eqno{(\EE.3)}
$$
Of course, it is enough to define the quadratic form
$$
B(v,v) \ \equiv \  {\cal B}(v-iJv, v-iJv)
\eqno{(\EE.3)'}
$$
which is real-valued and  determines (\EE.3) by polarization.  From (\EE.3)$'$ it is obvious that
$$
{\cal B}\ \geq\ 0 \qquad\iff\qquad B\ \geq\ 0.
\eqno{(\EE.4)}
$$

Now given any $B\in \Sym_\bbr(T)$, the {\bf hermitian symmetric part of $B$} is defined to 
be the element
$$
B^J(v,v) \ \equiv \  B(v,v) +B(Jv, Jv) 
\eqno{(\EE.5)}
$$
which belongs to $\HSymr T$. (Usually one inserts a $\half$ in (\EE.5), but here it is cleaner not to do so.)

Combining this algebra with (\BB.2) and adding $J$ to the notation, we have three isomorphic 
vector spaces
$$
\L^{1,1}(T(J)) \ \cong\  \HB(T_{1,0}(J)) \ \cong\ \HSymr {T(J)}. 
\eqno{(\EE.6)}
$$
It is important that the last vector space $\HSymr {T(J)}$ is a 
real vector subspace of $\Sym_\bbr(T(J))$.  We have denoted a triple of elements that 
corresond under these isomorphisms by $\b$, $\cb$ and $B$ respectively.  Note that
$\b\geq0 \ \iff\ \cb\geq0 \ \iff\ B\geq0$.  

We now  apply this to the case of the complex hessian of a function $\vf$ at a point $x\in X$
where $\b = i\partial\dbar \vf$ and $\cb = \ch(\vf)$.  The third element
$B\equiv H(\vf) \in\HSymr {T_x(J_x)}$ is called the {\bf real form} of the complex hessian.
Note that 
$$
 i\partial\dbar \vf\ \geq\ 0 \IFF \ch(\vf) \ \geq\ 0\IFF H(\vf) \ \geq\ 0.
\eqno{(\EE.7)}
$$

The formula (\BB.4)  provides a formula for $H(\vf)$.

\Lemma{\EE.1}  {\sl
The real form $H(\vf)$ of the complex hessian $\ch(\vf)$ is given by the polarization 
of the real quadratic form
$$
H(\vf) (v,v) \ =\bigl\{vv +(Jv)(Jv) + J\left ([v, Jv]\right)\bigr\} \vf
\eqno{(\EE.8)}
$$
defined for all real vector fields $v$ (where the vector fields act on functions in the standard way).
}

\pf As an operator on $\vf$,  $H(v,v)$ is given by 
$$
H(v,v) \ =\ \ch(v-iJv, v-iJv),
$$
which can be expanded out using  (\BB.4)$'$   to yield (\EE.8).\qed

\medskip

In euclidean space $\bbr^N$ with coordinates $t=(t_1,...,t_N)$,
let $p= D\vf$ (evaluated at $t$), and  let $D^2\vf$ denote both the second derivative matrix
$A\equiv (\!( {\partial^2 \vf  \over \partial t_i\partial t_j} )\!)$
as well as the quadratic form
$$
A(v,v) \ =\ (D^2\vf)(v,v) \ =\ \sum_{i,j=1}^N  {\partial^2 \vf  \over \partial t_i\partial t_j}(t) \,v_i v_j \qquad
{\rm where\ \ \  } v \ \equiv \ \sum_{j=1}^N  v_j {\partial   \over \partial t_j}.
$$
A  calculation gives the following.
\Prop{\EE.2}  {\sl
Suppose that $J$ is an almost complex structure on an open subset $X\ss \bbr^{2n}$.  
Let $v$ be a constant coefficient vector field on $X$  (i.e., $v= \sum_j v_j {\partial\over \partial t_j}$ where
the $v_j$'s are constants).  Then
$$
H(\vf)(v,v)\ =\ (D^2\vf)(v,v) + (D^2\vf)(Jv, Jv) + (D\vf)\bigl \{ (\n_{Jv}J)(v)  - (\n_vJ)Jv\bigr\}  
\eqno{(\EE.9)}
$$
(where $\n_vJ$ denotes the standard directional derivative of the matrix-valued function $J$).
Equivalently,
$$
H(\vf)\ =\  (A+E(p))^J  
\eqno{(\EE.9)'}
$$
where $p\equiv D\vf$, $A\equiv D^2\vf$ and the section $E\in \G(X, \Hom_\bbr(\bbc^n, \Symr(\bbc^n))$
is defined by
$$
E(p)(v,v) \ \equiv\  \bra {(\n_{Jv} J)(v)}  {p}
\eqno{(\EE.10)}
$$
where $\bra \cdot\cdot$ is the standard real  inner product on $\bbc^n$.}

\pf
By (\EE.8) 
$$
(H\vf)(v,v)\ =\ (D^2\vf)(v,v) + (D^2\vf)(Jv, Jv) + (D\vf)\bigl \{ (\n_{Jv} J)(v) +J[v,Jv]\bigr\}. 
$$
Now $[v,Jv] =  (\n_v J)v$, and $J^2=-I$ implies that $(\n_v J)J+J(\n_v J)=0$. Hence we have
$J[v,Jv] = J(\n_v J)v = -(\n_v J)Jv$, which proves (\EE.9).
To prove (\EE.9)$'$, 
note  that $A^J(v,v) =  A(v,v) +  A(Jv,Jv)$ and that $E(p)^J =  \bra {(\n_{Jv} J)(v)}  {p} - \bra {(\n_v J)(Jv)}  {p}$.\qed

 \medskip
 \noindent
 {\bf Remark.}
 Using the conventions $v=\sum v_je_j$, $Jv=\sum v_jJ_{jk}e_k$
 where $e_j= {\partial\over \partial t_j}$, the reader may wish to derive (\EE.9) and (\EE.9)$'$ from
 (\EE.8) using matrices.\medskip

We now continue with the  linear algebra.
Consider  a finite dimensional real vector space $T$ of dimension
$2n$ with two  almost complex structures $J$ and $J_0$,   
inducing the same orientation on $X$.
%(and for later use two volume forms $\l$ and $\l_0$,
%all of which are orientation compatible).  
Think of   $J_0$ as ``standard'' or  ``background'' data.
Assume we are given $g\in {\rm GL}^+(T)$ such that
$$
J \ =\  g J_0 g^{-1}
$$
Extend $g$ as usual to $GL_\bbc(T\otimes \bbc)$. The induced action on $\Sym_\bbr(T)$ 
is given by 
$$
(g^*B)(v,w) = B(gv,gw)
\eqno{(\EE.11)}
$$

The next result is more than we need but it should clarify the algebra.

\Lemma{\EE.3} {\sl
Consider elements 
$$
\a\in \L^{1,1}_\bbr (T, J),\quad \ch \in  \HB(T_{1,0}(J)),\quad {\rm and}\ \ 
 H\in \HB_\bbr(T, J) \ss \Sym_\bbr(T). 
$$
Then the pull-backs by $g$ satisfy
$$
g^*\a\in \L^{1,1}_\bbr (T, J_0),\quad  g^* \ch \in  \HB(T_{1,0}(J_0)),\quad {\rm and}\ \ 
 g^* H\in \HB_\bbr(T, J_0) \ss \Sym_\bbr(T). 
$$
Moreover, if the three elements $\a, \ch$ and $H$ correspond to each other under the isomorphisms
(\EE.6), then so do the elements $g^*\a, g^*\ch$ and $g^*H$. 
Any one of the six is $\geq0$ if and only if all six are $\geq0$. 
Finally, for any  $B\in \Sym_\bbr(T)$
%  the $J_0$-hermitian part
%of $g^* B$ equals $g^*$ of the $J$-hermitian symmetric part of $B$. That is,
$$
(g^* B)^{J_0}\ =\ g^*\left( B^J\right).
\eqno{(\EE.12)}
$$
}

\pf 
Since $Jg = gJ_0$, we have $g^* J^*= J_0^* g^*$. 
Finally, to prove (\EE.12) note that  $(g^*B)^{J_0} =
g^* B + J_0^* g^* B = g^* B + g^* J^*B = g^* (B+J^*B) = g^*(B^J)$ using (\EE.5).\qed
\medskip

Applying this to the complex hessian $\ch(\vf)$
and   its real form $H(\vf)$, yields
$$
\ch(\vf) \ \geq\ 0 \qquad\iff\qquad 
\left ( g^* H(\vf)\right)^{J_0} \ \geq\ 0.
\eqno{(\EE.13)}
$$

\Ex{\EE.4. (The Standard Complex Structure on $\bbc^n$)}
Let $J_0$ denote the standard complex structure ``$i$'' on $\bbc^n$.  With 
$V\equiv \sum_{j=1}^n c_j {\partial\over \partial z_j}$ the complex hessian $\ch_0$ is
  given by  $\ch_0(\vf)(V,V)  =
 \sum_{j,k=1}^n \left( {\partial^2 \vf\over  \partial z_j  \partial\bar z_k} \right)c_j \overline{c}_k$.
The real form $H_0$ of this complex hessian can be most succinctly expressed as
$$
H_0 \vf \ =\ (D^2\vf)^{J_0} \ =\    D^2\vf + J_0^* D^2\vf
\eqno{(\EE.14)}
$$
since $\nabla J_0\equiv 0$. That is, $H_0\bf$ is simply the $J_0$-hermitian symmetric part of 
the second derivative $D^2\vf$.
\medskip

Now the subequation $F(J_0)$ on an open subset $X$ of $\bbc^n$ is easily computed to be
$$
\eqalign
{
F(J_0)\ &=\ X\times \bbc^n\times \cp^\bbc\qquad\qquad{\rm where}\cr
\cp^\bbc \ \equiv \{ A &\in \Sym_\bbr(\bbc^n) :  A^{J_0} = A+J_0^*A \geq 0\}.
}
$$
By this definition of $\cp^\bbc$ the equivalence (\EE.13) can be rewritten as 
$$
\ch(\vf) \ \geq\ 0
\qquad\iff\qquad
g^*H(\vf) \ \in\ \cp^\bbc.
\eqno{(\EE.13)'}
$$

Sometimes it is convenient to refer to the subequation $F(J_0)$ simply as  $\cp^\bbc$
since
$$
\vf\ \ {\rm is\ } J_0\ {\rm subharmonic\ on\ \ } X \qquad\iff\qquad 
(D^2\vf)(x) \in \cp^\bbc\ \ \ \forall\, x\in X.
\eqno{(\EE.15)}
$$
\medskip

We now assume our   euclidean space $\bbr^{2n}$ to be 
equipped with a standard complex structure $J_0$ and write 
$\bbc^n = (\bbr^{2n}, J_0)$ as above.  We further assume that our 
variable almost complex structure $J$ can be written in the form
$J=g J_0 g^{-1}$ for a smooth map $g:X\to \GL^+_\bbr(\bbr^{2n})$.
(This can be arranged in a neighborhood of any point $x\in X$ by
choosing $J_0=J_{x}$.)  The next result describes $F(J)$ in  these
coordinates as a perturbation  of the standard subequation $\cp^\bbc$.
It is the key to the two main theorems in this paper.

\Prop{\EE.5}  {\sl
Suppose $g:X\to \GL_\bbr^+(\bbr^{2n})$ defines an almost complex structure 
$J\equiv g J_0 g^{-1}$ on an open subset $X\ss \bbc^n$.  Let $\ch$ denote the
complex hessian for $J$.  Then 
$$
\ch(\vf) \geq 0 \qquad  \iff \qquad  g^* D^2\vf + g^*E(D\vf) \ \in\ \cp^\bbc
\eqno{(\EE.16)}
$$
where $E$ is defined by (\EE.10).
}

\pf
By (\EE.9)$'$, $H(\vf)$ is the $J$-hermitian part of $D^2\vf + E(D\vf)$.
Hence, by Lemma \EE.3, $g^*H(\vf)$ is the $J_0$-hermitian part of 
$g^*(D^2\vf + E(D\vf))$.  Therefore,
$$
g^*H(\vf) \in \cp^\bbc \qquad  \iff\ \qquad 
 g^* \left (D^2\vf +E(D\vf) \right)  \in \cp^\bbc.
$$
Combining this with (\EE.13)$'$ completes the proof of (\EE.16).\qed

%%%%%%%%%%%%%%%%%%%%%%%%%%%%%%%%%%%%%%%%%%%%%%%
%%%%%%%%%%%%%%%%%%%%%%%%%%%%%%%%%%%%%%%%%%%%%%%
%%%%%%%%%%%%%%%%%%%%%%%%%%%%%%%%%%%%%%%%%%%%%%%
%%%%%%%%%%%%%%%%%%%%%%%%%%%%%%%%%%%%%%%%%%%%%%%
%%%%%%%%%%%%%%%%%%%%%%%%%%%%%%%%%%%%%%%%%%%%%%%

%\vfill\eject
\vskip .3in

\noindent{\headfont  \FF. Restriction of \FPSH Functions.}\medskip

In this section we prove (in Theorem \FF.2) that the restriction of a \FPSH 
function to an almost complex submanifold is also plurisubharmonic
(as a function on the submanifold).  The difficulty of this result is somewhat surprising.
First we establish some easier facts.
\vfill\eject

\Prop{\FF.1} {\sl 
Suppose  $\Phi :(X',J') \to (X,J)$ is a pseudo-holomorphic map. 
\medskip\noindent
{(1)}  \ \ If $u\in C^2(X)$ is \FPSH on $(X,J)$,
 then $u\circ \Phi$ is \FPSHP{J'} on $(X',J')$.

\medskip\noindent
{(2)} \ \ If $j\in \FxJ$, then $\Phi^*(j) \in  \FJXP{x'}{J'}$
where $x=\Phi(x')$, that is, 
$$
\Phi^*  \left( \FJ \right) \ \ss\  \FJP{J'},
$$
}

\pf Since $\Phi^*$ commutes with $i\partial\dbar$, 
the pull-back of a smooth \FPSH function is $F(J')$-plurisubharmonic.  The proof of (2) is similar.\qed

\medskip

Now we state the more difficult result.

\Theorem{\FF.2. (Restriction)} {\sl
Suppose that $(X',J')$ is an almost complex submanifold of $(X,J)$.
If $u\in\USC(X)$ is \FPSH on $X$, then $u\bigr|_{X'}$ is 
\FPSHP{J'} on $X'$.}
\medskip

\pf  Since the result is local, we may choose coordinates which reduce us to the following
situation.  Suppose that $J$  is an almost complex structure on a neighborhood $X$ of the 
origin in $\bbc^n$, which agrees with the standard complex structure $J_0$ at $z=0$.
Suppose further that  $X'=(\bbc^m\times \{0\})\cap X$ is a $J$-almost-complex submanifold.
By shrinking $X$ if necessary we can find a smooth
mapping  $g: X \to \GL^+_{\bbr}(\bbc^n)$ with 
$$
g(0)={\rm Id}
 \and
J=g J_0 g^{-1}\quad {\rm  on\ }X.
\eqno{(\FF.1)}
$$    
Block the transformation $g$   as
$$
g\ =\ \left ( \matrix{   g_{11}    &    g_{12}   \cr    g_{21}    & g_{22}       \cr}\right)
\qquad
{\rm with \ respect\  to }\ \ \ \ \bbc^n=\bbc^m\times \bbc^{n-m}.
$$
Because of the next result,
we can choose $g$ such that 
$$
g_{21}\ \equiv \ 0   \quad{\rm on\ the\ submanifold\ } X'.
\eqno{(\FF.2)}
$$

\Lemma{\FF.3} {\sl
By further shrinking   $X$ if necessary, the mapping $g$ can be chosen to be of the form 
$$
g\ =\ {\rm I}+f \quad {\rm with\ } f \ {\rm complex\ antilinear.}
$$
With this choice, $f$ is unique and $f_{21} (\equiv g_{21}) \equiv 0$ on $X'$.
}

\def\I{{\rm I}}

\pf Recall that each $g\in \End_\bbr(\bbc^n)$ has a unique decomposition $g=h+f_1$ with
$h\in \End_\bbc(\bbc^n)$ complex linear, and 
$f_1\in \overline{\End}_\bbc(\bbc^n)$ complex anti-linear.
Since $h(0)=  \I$ , we may assume, by shrinking $X$,  that $h(x)$ is invertible for each $x\in X$.
Define $f\equiv f_1 h^{-1}$.  Then (since $h$ and $J_0$ commute), we have
$$
J \ =\  gJ_0 g^{-1} \  = \    gh^{-1}J_0 hg^{-1}     \ =\ 
g h^{-1} J_0 ( g h^{-1})^{-1}
\and
g h^{-1} \ =\ {\rm I}+f.
$$
This proves the first assertion.

For the uniqueness statement, suppose that $J = (\I+f_1)J_0(\I+f_1)^{-1}
=  (\I+f_2)J_0(\I+f_2)^{-1}$ with both $f_1$ and $f_2$ complex anti-linear.
Then $(\I+f_2)^{-1}(\I+f_1)$ commutes with $J_0$, i.e., is complex linear.
However, the complex anti-linear part of 
$
(\I+f_2)^{-1}(\I+f_1) = (\I+f_2^2)^{-1}(\I-f_2)(\I+f_1)
$
is 
$
 (\I+f_2^2)^{-1}(f_1-f_2),
$
and so $f_1-f_2\equiv0$. 

It remains to prove  the last assertion. To begin we block $J$ with respect to the splitting
$\bbc^n=\bbc^m\times \bbc^{n-m}$ as above.  Then since $X'$ is an almost complex
submanifold, the component $J_{21}$ must vanish along $X'$.  Therefore, the 21-component
of $Jg$ equals $J_{22}g_{21} = J_{22}f_{21}$, while the 21-component of $gJ_0$ equals
$g_{21}i=f_{21}i = -i f_{21}$.  Since $Jg=g J_0$, this proves that $(J_{22} +i) f_{21}=0$.
Finally, $J_{22}(0)=i$, so that $f_{21}=0$ along $X'$ near the origin.\qed

\medskip
\noindent
{\bf Note.}  The last two statements can be seen in another way.
An almost complex structure $J_0$ on a real vector space $T$ is equivalent to
a decomposition $T\otimes_\bbr \bbc = T_{1,0}\oplus T_{0,1}$ with 
$ T_{0,1}  =  \overline{ T_{1,0}}$.  Another complex structure $J$
on $T$, inducing the same orientation, has a similar decomposition,
and $T_{1,0}(J)$ is the graph in $T_{1,0} \oplus T_{0,1}$ 
of a unique complex linear map $f:T_{1,0}\to T_{0,1}$
(or equivalently, a $J_0$-anti-linear map $f:T\to T$).
Suppose now that  $S\ss T$ is  a $J$-complex subspace which is also $J_0$-invariant.
Then ${S_{1,0}}(J) = S_{1,0}$  and so $f\bigr|_{S_{1,0}} = 0$.
This is  the condition that $f_{21}=0$.

%Equivalently,
%$T_{1,0}(J)$ is the image  in $(T,J_0)\oplus (T,-J_0)$ of a unique map
%of the form $\I\oplus f$ where $f:T\to T$ and $f \circ J_0= -J_0\circ f$. 
%With $J_0$   the standard complex structure, 
%this mapping $f$ is the one constructed above.

\medskip
\noindent
{\bf Completion of the Proof of Theorem \FF.2.} 
We apply Proposition \EE.5.
 First note that 
(\EE.16) can be restated, using $p=Du$ and $A=D^2u$ as 
$$
(p,A) \in F_x(J) \qquad \iff \qquad g(x)^t (A  + E_x(p)) g(x) \in \cp^\bbc(\bbc^n)
\eqno{(\EE.16)'}
$$
Set $L_x(p)\equiv g^t(x) E_x(p) g(x)$.

\Lemma {\FF.4}  {\sl
If $p=(p',p'') \in \bbc^n=\bbc^m\times \bbc^{n-m}$, 
then for $x\in  X' \equiv \bbc^m\times \{0\}$,   $L_x((0,p'')) \in \Sym_\bbr(\bbc^n)$
 vanishes when restricted to $X'$
as a quadratic form.
}
\pf  
Suppose that $v$ is a vector field tangent to $X'$
along $X'$. Then since $X'$ is an
almost complex submanifold, the vector field $Jv$ also has this
property.  It now follows directly from  (\EE.10) that if  $p=(0,p'')$, then $E(p)(v,v)\equiv0$
along $X'$.  In other words the component 
$$
E_{11}(0,p'') \ \equiv\ 0 \qquad{\rm along}\ \  X'.
$$
This together with (\FF.2) implies that 
$$
\qquad\qquad\   L_{11}(0,p'') \ \equiv\ 0 \qquad{\rm along}\ \  X'. \qquad\qquad\mathqed
$$

\medskip

The hypotheses (7.12) of the 
 Restriction Theorem 8.1 in [HL$_3$] are now established.   
To check this, note  first that the matrix function $h$ in (7.12) of [HL$_3$]  equals our $g^t$.
Hence, by  (\FF.2) we have $h_{12} = g_{21} =0$ on $X'$ as required.
The $g$ in (7.12) of [HL$_3$]  is taken to be the identity, so its (12)-component vanishes on $X'$.
 The last  part of  (7.12) in [HL$_3$]  is exactly   Lemma \FF.4 above.
Theorem \FF.2  now follows from
  Theorem 8.1 in [HL$_3$]. \qed
\medskip

     For 1-dimensional almost complex manifolds  Theorem \FF.2 has a converse,
     which we investigate in the next section.

%%%%%%%%%%%%%%%%%%%%%%%%%%%%%%%%%%%%%%%%%%%%%%%
%%%%%%%%%%%%%%%%%%%%%%%%%%%%%%%%%%%%%%%%%%%%%%%
%%%%%%%%%%%%%%%%%%%%%%%%%%%%%%%%%%%%%%%%%%%%%%%
%%%%%%%%%%%%%%%%%%%%%%%%%%%%%%%%%%%%%%%%%%%%%%%
%%%%%%%%%%%%%%%%%%%%%%%%%%%%%%%%%%%%%%%%%%%%%%%

\vfill\eject
%\vskip .3in

\noindent{\headfont  \II. The Equivalence of \FPSH  Functions 
and  Standard Plurisubharmonic Functions.}\medskip

  In complex dimension one, each almost complex structure is integrable, i.e., each almost complex manifold $(\Sigma, J)$ of real dimension 2 is a Riemann surface. There are many equivalent 
definitions for subharmonic functions on a Riemann surface, and  Definition \CC.3 is one of these.
We assume these facts without further discussion.
 
 A ``standard''  definition of a plurisubharmonic function on a 
complex manifold makes perfect sense on an almost complex manifold.

\Def{\II.1}  An upper semi-continuous function $u$ on an almost complex manifold
$(X,J)$ is said to be {\bf   plurisubharmonic in the standard sense}  if  its restriction  to each
holomorphic curve $\Sigma$ in $X$ is subharmonic.
     \medskip
     
   The Restriction  Theorem \FF.2  implies the forward implication in the next 
   result.  The abundance of holomorphic curves on an almost complex manifold will be used to prove the reverse implication.
     
     \Theorem{\II.2}  Given $u\in \USC(X,J)$ on an almost complex manifold $(X,J)$, 
     $$
     u\ \ {\rm is \ } F(J)\ {\rm plurisubharmonic} \qquad \iff\qquad 
      u\ \ {\rm is \    plurisubharmonic \ in\ the\ standard\ sense}.
     $$     
\medskip

Consequently, we may simply call these functions {\bf $J$-plurisubharmonic}.
The set of all such functions on $X$ will be denoted by $\PSH^J(X)$.

In this section we shall replace the notation $\FJ$ by $F^X$ to emphasize the manifold and to
suppress confusion with notation for 2-jets.  The jet version of Theorem \II.2 can be stated as follows.

\Lemma{\II.3}  
$$
\bbj \in F^X_z \quad\iff\quad i^*(\bbj) \in F^{\Sigma}_\z \quad{\rm for\ all\ holomorphic \ curves\  } i:\Sigma  \to X\ 
{\rm with\ } i(\z)=z.
$$

 \pf ($\Rightarrow$): This is the special case of Proposition \FF.1(2) saying that 
 $i^*(F_z^X)\ss F_\z^\Sigma$.

  ($\Leftarrow$):
 Pick a smooth function  $\vf$  with $J_z(\vf) \equiv \bbj$.  Assume that
 $ i_{\Sigma}^*(\bbj)  \in  F^{\Sigma}_\z$ for all $i_{\Sigma}:\Sigma\to X$ and $i_{\Sigma}(\z)=z$.
     We must show that $\ch_z(V,V)(\vf)\geq 0$ for all $V\in (T_{1,0}X)_z$.  By [NW] there exists
     $i:\Sigma \to X$ with $i(\z)=z$ and $i_*({\partial\over \partial \z})=V$.
     Hence, $\ch_z(V,V)(\vf) = \ch_z(i_*({\partial\over \partial \z}),i_*({\partial\over \partial \z}))(\vf)
     = \ch_{\z}({\partial\over \partial \z}, {\partial\over \partial \z})(\vf \circ i)$ which is $\geq 0$ since
     $J_{\z}(\vf\circ i) = i^*(J_z\vf) \in F^\Sigma_\z$.\qed

\medskip
\noindent
{\bf Proof of Theorem \II.2}  As noted above we only need to prove $\Leftarrow$.
Assume that  $u\in\USC(X)$ and that 
$u\circ i$ is subharmonic on $\Sigma$ for each
holomorphic curve $i:\Sigma \to X$.  Pick a point $z\in X$ and a test function $\vf$ 
     for $u$ at $z$.  Suppose $i:\Sigma \to X$ is a pseudo-holomorphic curve 
   with $i(\z)=z$. Obviously, $\vf\circ i$ is a test function for $u\circ i$ at $\z$. Since $u\circ i$ is
   subharmonic, $J_{\z}(\vf\circ i)\in F^\Sigma_\z$.  By Lemma \II.3 this is enough to imply that
   $J_z(u)\in F_z^X$, since $i^*(J_z(\vf))= J_{\z}(\vf\circ i)$.\qed

       \medskip
     
     There are advantages to using the concept of \FPSH over the standard one
     on an almost complex manifold.  This is apparent for example, in the 
     Section \HH\  on the Dirichlet Problem.
     
     We conclude the section by mentioning a more elementary application illustrating
     the abundance of plurisubharmonic functions locally. We say that a system of local
     coordinates  $z=(z_1,...,z_n)$ on $X$ is {\sl standard at} $x\in X$ if $z(x)=0$ and 
     $J_x \cong J_0 (=i)$.

\Prop{\II.4} {\sl
Suppose $z$ is a local coordinate system for $(X,J)$ which  is standard at $x$.
Then the function $u(z)=|z|^2$ is strictly \FPSH on a neighborhood of $x$.
}
   \pf
    By (\EE.9), we see that, since $Du_0=0$, the real form of the $J$-hermitian hessian
at the origin is $H(u)_0 = 2I$.
   Hence, $H(u)$ is positive definite in a neighborhood of 0.   \qed
   
   \Cor{\II.5} {\sl
   Each point of an almost complex manifold has a neighborhood system of domains with strictly
   pseudo-convex smooth boundaries.
   }
     \medskip
     More precisely, these boundaries are strictly $\FJ$-convex in the sense of [HL$_2$].

%%%%%%%%%%%%%%%%%%%%%%%%%%%%%%%%%%%%%%%%%%%%%%%
%%%%%%%%%%%%%%%%%%%%%%%%%%%%%%%%%%%%%%%%%%%%%%%
%%%%%%%%%%%%%%%%%%%%%%%%%%%%%%%%%%%%%%%%%%%%%%%
%%%%%%%%%%%%%%%%%%%%%%%%%%%%%%%%%%%%%%%%%%%%%%%
%%%%%%%%%%%%%%%%%%%%%%%%%%%%%%%%%%%%%%%%%%%%%%%

%\vfill\eject
\vskip .3in

\noindent{\headfont  \HH. The Dirichlet Problem.}\medskip

In this section we consider the Dirichlet problem for the subequation $\FJ$  and for the more general ``inhomogeneous''  subequation $F(J,f)$ defined by adding the
condition $(i\partial\dbar u)^n \geq f \l$ where $\l$ is a fixed  volume form on $X$ and 
$f \in C(X)$  satisfies $f\geq0$.  Recall that by
 Definition \CC.1,    $F_x(J) = \{ J_x^2 u : i\partial\dbar u \geq 0 \ {\rm for}\ u\ {\rm smooth\ near }\ x\}$,
 which we state more succinctly by saying
  $$
 F(J) \ \ {\rm is\ defined\ by\ }\ \  i \partial\dbar u\ \geq\ 0.
 \eqno{(\HH.1)}
 $$
We now fix a (real) volume form $\l$ on $X$ compatible with the orientation induced by $J$. Then
each $f\in C(X)$ with $f\geq0$ determines a subequation 
$$
 F(J, f) \ \ {\rm  defined\ by\ }\ \ \   i \partial\dbar u\ \geq \ 0 \and (i\partial\dbar u)^n \ \geq \ f \l.
 \eqno{(\HH.2)}
 $$

For simplicity we abbreviate $F\equiv F(J,f)$.

\Def{\HH.1} A smooth function $u$ on $(X,J)$ is  {\bf $F$-harmonic} if $J^2_x(u) \in \partial F_x$
for each point $x\in X$.

\Ex{\HH.2. (The Standard Model)} 
In the model case $(X,J)=(\bbc^n , J_0)$ with $\l= \l_0 = 2^n \times$ the 
standard volume form on $\bbc^n$ and $f_0\geq 0$ a continuous function, we have
$$
D^2 u\ \in\   {\bf F}(J_0, f_0) 
\qquad\iff\qquad
i \partial\dbar u \geq 0 \ \ {\rm and} \ \ (i\partial\dbar u)^n  
 \  \geq f_0 \l_0,
 \eqno{(\HH.3)}
 $$
 where the second inequality can be rewritten as
 $$
  \det_\bbc\left({\partial^2 u\over \partial z\partial \overline z}\right) \ \geq\ f_0. 
 $$
Furthermore,  
$$
D^2 u\ \in\ \partial {\bf F}(J_0, f_0) 
\quad\iff\quad
i \partial\dbar u \geq 0 \ \ {\rm and} \ \ (i\partial\dbar u)^n 
 \  =\  f_0 \l_0,
 \eqno{(\HH.4)}
 $$
 where the  equality can be rewritten as
 $$
  \det_\bbc\left({\partial^2 u\over \partial z\partial \overline z}\right) \  = \ f_0.
 $$

 This model subequation $\bbf \equiv \bbf(J_0,f_0) \ss\Sym(\bbr^{2n})$
 is pure second-order, but not constant coefficient unless $f_0$ is constant.  Rewriting (\HH.3) 
 using jet coordinates we have that 
 $$
 A\in \bbf_x(J_0,f_0) \qquad\iff\qquad A_\bbc \geq 0 \quad{\rm and}\quad \det_\bbc \left(A_\bbc \right) \ \geq f_0(x)
 \eqno{(\HH.3)'}
 $$
\medskip
 
 The notion of $F$-harmonicity 
can be extended to $u\in \USC(X)$ as follows.  Define the {\bf Dirichlet dual} of 
$F$ to be the set
$$
\wt F \ =\ \sim\left(-\Int F  \right)\ =\  -\left(\sim\Int F \right).
\eqno{(\HH.5)}
$$
One can show that $\wt F$ is also a subequation, i.e., it is a closed subset of
 the (reduced) 2-jet bundle ${J}^2(X)$ satisfying conditions (P) and (T)
 (cf. [HL$_2$] and  Definition \CC.5 above).

 Note that 
  $$
F \cap \left( - \wt F \right) \ =\ F \cap  \left( \sim\Int   F \right) \ =\ 
\partial F.
\eqno{(\HH.6)}
$$   
For any subequation $F$ the  notion of an upper semi-continuous $F$-subharmonic function 
is defined by using test functions (see Definition \CC.6),
and Definition \HH.1 can be extended to continuous functions as follows.

 \Def{\HH.3} A  function $u \in C(X)$     is $F$-{\bf harmonic} if $u$ is  
$F$-subharmonic  and $-u$ is $\wt F$-subharmonic on $X$.
     \medskip
     
The dual $\wt F$ of $F\equiv F(J,f)$ is defined fibre-wise by:
$$
\eqalign
{{\rm either} \ \ &(1) \ \  i\partial \dbar(-u)\geq0 \quad
 {\rm and }\quad \left(  i\partial \dbar(-u)\right)^n \geq f(x) \l   \cr
\quad \  \  {\rm or} \ \ &(2)\ \   i\partial \dbar(-u) \not\geq0. \cr
}
\eqno{(\HH.7)}
$$

     Now suppose that $\O$ is an open set in $X$ with smooth boundary $\partial \O$
     and that $\overline \O=\O\cup \partial \O$ is compact.
     
     \Def{\HH.4} We say that {\bf uniqueness holds} for the $F$-Dirichlet problem (DP) on
     $\O$ if  given  $\vf\in C(\bo)$ and $v,w\in C(\ob)$ with  $v$ and $w$ both $F$-harmonic on $\O$,
     then   
     $$
     v\ =\ w\ =\ \vf \ \ {\rm on}\ \bo \qquad\Rightarrow\qquad v\ =\  w \ \ {\rm on}\ \O.
     $$
     
     We say that {\bf existence  holds} for   (DP) on
     $\O$ if  given  $\vf\in C(\bo)$, the Perron function 
    $$
    U\ \equiv\  \sup_{\cf(\vf)} u \qquad{\rm where} \qquad \cf(\vf) \ \equiv \ 
    \left\{u\in \USC(\ob) \cap F(\O) : u\bigr|_{\bo} \leq\vf \right\}
    $$ 
     satisfies
     \medskip
     
     (1)\ \  $U\in C(\ob)$, 
     \quad (2) \ \ $U$ is $F$-harmonic on $\O$ \ \ and\ \ 
     (3)\ \  $U = \vf$ on $\bo$.

     \medskip
     
     \Theorem{\HH.5. (The Dirichlet Problem)} {\sl 
     For the subequation    $F=F(J, f)$ on an almost complex manifold $(X,J)$, we have the following.

     Uniqueness holds for the $F$-Dirichlet problem on 
   $(\O,\bo)$  if the  almost complex manifold $(X,J)$  supports a $C^2$ 
   strictly $F(J)$-plurisubharmonic function.
     
     Existence holds for the $F$-Dirichlet problem if $(\O,\bo)$ has a strictly $F(J)$-plurisubharmonic defining function.
     }

     \Cor{\HH.6}  {\sl  On any almost complex manifold $(X,J)$ each point has a
      fundamental neighborhood system of domains $(\O,\bo)$ for which both
      existence and uniqueness hold for the Dirichlet problem.}
      
      \medskip
      \noindent
      {\bf Proof.} The  idea is to apply the results of [HL$_2$]  to the subequations $F =F(J,f)$ despite the lack of a riemannian metric.  Proposition \EE.5 above states that $F(J)$
      is locally   jet equivalent  to $\cp^\bbc = {\bf F}(J_0)$, the standard constant coefficient subequation
      on $\bbc^n$, while, more generally,  Proposition \HH.8  below states that $F(J, f)$ is locally jet-equivalent to 
      ${\bf F}(J_0, f_0)$, the subequation described in Example \HH.2.

\medskip
\noindent
{\bf Uniqueness.}   This is a standard consequence of comparison.
\Theorem {\HH.7. (Comparison)}   
  {\sl  If $(X,J)$ supports a $C^2$ strictly $J$-plurisubharmonic function, then
  comparison holds for $F\equiv F(J,f)$, that is, for all $u\in F(X)$, $v\in {\wt F}(X)$, and 
  compact subsets $K\ss X$: }
  $$
 u+v\ \leq \ 0 \ \ \ {\rm on}\ \partial K 
 \qquad\Rightarrow\qquad
 u+v\ \leq \ 0 \ \ \ {\rm on}\  K.
 $$ 
\pf
Section 10 in [HL$_2$] contains a proof that any subequation $F$ which is locally
      affinely jet equivalent to a constant coefficient subequation
     ${\bf F}$ must satisfy {\sl local weak comparison}.   
     Even though $\bbf(J_0,f_0)$ is not constant coefficient, it is sufficiently close that the 
     Theorem on Sums used in [HL$_2$] can be used to prove  local  weak comparison for 
     $F(J,f)$.  This is done in Proposition \HH.9 below.
     Theorem 8.3 in [HL$_2$] states that {\sl local 
      weak comparison}  implies {\sl weak comparison} on $X$.  Thus,
      $$
      {\rm Weak \ comparison \ holds\ for \ } F(J,f) \ {\rm on\ any \ almost\  complex\  manifold}\ (X,J).
      \eqno{(\HH.8)}
      $$

It is easy to see that under fibre-wise sum
$$
F(J,f) + F(J) \ \ss\ F(J,f),
 \eqno{(\HH.9)}
$$
that is, the convex cone subequation $F(J)$ is a monotonicity cone for $F(J,f)$.

      Theorem 9.5 in [HL$_2$] states that if there exists a $C^2$ strictly $F(J)$-plurisubharmonic function
      on $X$, then the {\sl strict approximation}  property holds for $F(J,f)$ on $X$.
      Finally, Theorem 9.2 in [HL$_2$] states that for any subequation $F$, if both weak comparison
      and strict approximation hold, then {\bf  comparison holds} on $X$.
Modulo proving Propositions \HH.8 and \HH.9 below, this completes the proof of 
comparison. \qed

\medskip      \noindent
{\bf Existence.}  This  is a consequence of Theorem 12.4 in [HL$_2$], since
 it is straightforward to see from the definition that strict boundary convexity for $F(J)$ is the 
same as strict boundary convexity for $F(J,\l)$.  \qed

 \medskip
 In the   case $F=F(J, f)$, if $w$ is $C^2$ and $(i\partial\dbar w)^n \leq  f \l$,
 i.e., $J^2 w\notin \Int F$, then $v\equiv -w \in\wt F(X)$ and comparison applies to $v,u$
for any $u\in F(X)$.  In other words, $u\leq w$ on $\partial K \ \Rightarrow\ u\leq w$ on $K$.
Therefore, as a special case of comparison we have the following.
$$
\eqalign
{
&{\rm For\ all\ }\ u\ {\rm which\ are\ } J\,{\rm psh \ and\  satisfy \ } 
(i\partial\dbar u)^n \geq f \l\ {\rm in \ the\  viscosity\ sense,}    \cr
&{\rm if \ } w\in C^2(X) \ {\rm satisfies}\ 
(i\partial \dbar w)^n \leq f \l \ {\rm and \ } u\leq w \ {\rm on}\ \partial K,  
 \ {\rm then\ } u\leq w \ {\rm on}\ K.
 }
\eqno{(\HH.10)}
$$

\medskip
\centerline{\bf Local Jet-Equivalence with the Standard Model.}
\medskip

Recall the standard model in Example \HH.2.
 
\Prop{\HH.8} 
{\sl
Choose local coordinates in $\bbr^{2n} \cong \bbc^n$ as in Section 4 above.
In these coordinates  the
 subequation $F(J, f)$ is locally jet-equivalent to the standard model
subequation ${\bf F}(J_0, f_0)$ defined by 
$$
\det_\bbc\left( {\partial^2 u\over \partial z\partial \overline z}  \right) \ \geq\ f_0
\and
  {\partial^2 u\over \partial z\partial \overline z}  \ \geq\ 0,
$$
Here $f_0 = \b f$ where $\b>0$ is a smooth function, independent of $f$.}

\medskip

By the definition of jet-equivalence this Proposition states that there exists a 
GL$^+(\bbr^{2n})$-valued smooth function $h$ and a $(\bbr^{2n})^*\times \Sym(\bbr^{2n})$-valued
smooth function $L$ such that, with jet coordinates $p=Du$ and $A=D^2u$,
$$
(p,A) \in F_x(J,f) 
\qquad\iff\qquad
h(x) A h^t(x) + L_x(p) \ \in\ \bbf_x(J_0, f_0).
\eqno{(\HH.11)}
$$

\pf   
We recall the  linear algebra from \S \EE \ (cf. Lemma \EE.3) which involves
 a finite dimensional real vector space $T$ of dimension
$2n$ with two  almost complex structures $J$ and $J_0$,   and two volume forms $\l$ and $\l_0$,
all of which are orientation compatible.  
(We think of   $J_0$ and $\l_0$ as ``standard'' or  ``background'' data.)
Assume we are given $g\in {\rm GL}^+(T)$ such that
$$
J \ =\  g J_0 g^{-1}
$$
Let $\cp^\bbc(T,J) \ss \Sym_\bbr(T)$ consist of all $H\in \Sym_\bbr(T)$ such that:
\medskip
\centerline{
(1)\ $H$ is $J$-hermitian, i.e., $J^* H = H$ \and
(2) \ $H\geq0$.}
\medskip

\noindent
Technically speaking, $\cp^\bbc(T,J)$ is not a subequation, but it can be used to define
$F(J)$ by setting
$$
F_x(J)\ \equiv\ \{J^2_x \vf : H_x(\vf) \in\  \cp^\bbc(T,J)\}.
$$
Using this notion we have, as a restatement of Proposition \EE.5, 
where $A\equiv D_x^2\vf$ and $p\equiv D_x\vf$, that 
$$
(p,A) \in F_x(J)
\qquad\iff\qquad
g(x)\left(  A + J_0^* A  \right) g^t(x) + L_x(p) \ \in\ \cp^\bbc(\bbc^n, J_0)
\eqno{(\HH.12)}
$$
where $L_x(p) \equiv g(x)(E_x(p) +J_0^* E_x (p)) g^t(x)$ with $E$ defined by (\EE.8).
This is the statement that  $F(J)$ and $F(J_0)$ are  jet-equivalent.

Now with the notation of  Lemma \EE.3, set $\a=i\partial\dbar \vf$, $\ch\equiv \ch(\vf)$
and $H\equiv H(\vf)$. Then $H'\equiv g^* H = g(x)\left(  A + J_0^* A  \right) g^t(x) + L_x(p)$
is the expression in (\HH.12), and (\HH.12) is the statement that 
$H'\geq0$ and $H'$ is $J_0$-hermitian.

Let $\b>0$ be defined by
$$
g^* \l \ =\ \b \l_0.
\eqno{(\HH.13)}
$$
Then $g^*(f\l) = f\b \l_0 = f_0\l_0$ where $f_0\equiv \b f$.  Therefore,
$$
\a^n \geq f \l 
\qquad\iff\qquad
(g^*\a)^n \geq f_0 \l_0
\qquad\iff\qquad
\det_\bbc A_{\bbc}' \geq f_0
\eqno{(\HH.14)}
$$
where $A_{\bbc}'  \in M_n(\bbc)$ is the matrix representative of $\a' = g^*\a \in \L^{1,1}(T,J_0)$
with respect to any volume-compatible (i.e., $\l_0$-compatible) basis of $T_{1,0}(J_0)$.
This completes the proof of Proposition \HH.8.\qed

\vskip.3in
\centerline
{\bf  Local Weak Comparison for the Standard Model}
\medskip

The Theorem on Sums can be used to prove a weaker form of comparison
for $\bbf \equiv   \bbf(J_0,f)$ defined by:
$$
A\in \bbf_x
\IFF
A_\bbc \ \geq \ 0\and \det_\bbc A_\bbc \ \geq \ f(x),
\eqno{(\HH.15)}
$$
as well as for $F \equiv F(J, f)$ defined by:
$$
(p,A) \in F_x
\IFF
g(x) A g^t(x) + L_x(p) \in \bbf_x.
\eqno{(\HH.16)}
$$

A notion of strictness which is uniform is employed.  Given $c>0$, 
define $\bbf^c$ by 
$$
A\in \bbf_x^c
\IFF
B(A; c) \ \ss\ \bbf_x
\eqno{(\HH.17)}
$$
where $B(A; c)$ denotes the ball  in $\Sym_\bbr(\bbc^n)$ about
$A$ of radius $c$.  Define $F^c$ by 
$$
(p,A) \in F_x^c
\IFF
g(x) A g^t(x) + L_x(p) \in \bbf_x^c.
\eqno{(\HH.18)}
$$
Thus $F$ is jet-equivalent to $\bbf$,  and $F^c$ is jet-equivalent to $\bbf^c$.
Fix an open set $U\ss \bbc^n$.

\Prop{\HH.9} {\sl
If $u\in F^c(U)$ and $v\in \wt F(U)$, then $u+v$ satisfies the maximum principle.
}
\pf 
If the maximum principle fails for $u+v$, then (by [HL$_2$], Theorem C.1) 
the following quantities exist and have the stated properties:
$$
\eqalign
{
&(1) \ \  z_\e\ =\ (x_\e, y_\e) \ \to\ (x_0, x_0),  \cr
&(2) \ \  (p_\e, A_\e) \in F_{x_\e}^c \and (q_\e, B_\e) \in { \wt F}_{y_\e},     \cr
&(3) \ \  p_\e\ =\ {x_\e - y_\e  \over \e} \ =\ -q_\e  \and  { \ | x_\e - y_\e\ |^2 \over \e} \ \to\ 0, \cr
&(4) \ \   -{3\over \e} I \ \leq\    \left( \matrix{A_\e & 0 \cr 0 & B_\e} \right) 
\ \leq\   
{3\over \e}   \left( \matrix{I  & -I \cr -I & I} \right).    \cr
}
$$

It is straightforward to show that
$$
(q,B) \in {\wt F}_y \IFF g(y) B g^t(y) +L_y(q) \in {\wt \bbf}_y.
\eqno{(\HH.19)}
$$
(More generally, see Lemma 6.14 in [HL$_2$].) Thus (2) can be written as 
$$
\eqalign
{
(2)'\ \ &A_{\e}' \ \equiv \ g(x_\e) A_\e g^t(x_\e) + L_{x_\e}(p_\e)\ \in\ \bbf_{x_\e}^c, \ \ {\rm and}\qquad \ \ \ \  \cr
&B_{\e}' \ \equiv \ g(y_\e) B_\e g^t(y_\e) + L_{y_\e}(q_\e)\ \in\ {\wt \bbf}_{y_\e}.
}
$$
By definition of the dual, we have  $B_{\e}'  \in {\wt \bbf}_{y_\e}$ if  $-B_{\e}'  \notin \Int  {\bbf}_{y_\e}$, that is:
$$
\eqalign
{
(2)''\ \   {\rm either}\ \  &(i)\ \ -B_{\e}'  \notin \cp^\bbc \cr
 {\rm or}\ \  &(ii)\ \ -B_{\e}'  \in \cp^\bbc \ \ {\rm but}\ \ \det_\bbc\left( -B_{\e}' \right)_\bbc\ \leq \ f(y_\e).
}
$$
Recall that $\cp^\bbc \equiv \{B :  B_\bbc \geq 0\}$.

The calculation on page 442 of [HL$_2$] proves that there exist $P_\e\geq0$ and a number
$\L\geq0$ such that 
$$
 g(x_\e) A_\e g^t(x_\e) + g(y_\e) B_\e g^t(y_\e)  + P_\e \ =\ {\L\over \e} |x_\e-y_\e|^2.
\eqno{(\HH.20)}
$$
Setting ${A_\e}'' \equiv A_\e' +P_\e$, this proves
$$
(5)\ \ -B_\e' \ =\ {A_\e}'' - {\L\over \e} |x_\e-y_\e|^2 - L_{x_\e}(p_\e) -L_{y_\e}(q_\e).
$$
By positivity, $A_\e' \in \bbf_{x_\e}^c \ \Rightarrow\ {A_\e}''  \in \bbf_{x_\e}^c$, so that
$$
(6) \ \ \dist\left({A_\e}'', \sim \bbf_{x_\e}     \right) \ \geq\ c. \qquad\qquad\qquad \qquad\qquad \quad 
$$
This implies
$$
\eqalign
{
(6)' \ \ & (i) \ \ \dist\left({A_\e}'', \sim \cp^\bbc     \right) \ \geq\ c \ \ {\rm and} \qquad\qquad\qquad  \ \ 
\cr
& (ii) \ \ \det_\bbc\left({A_\e}'' \right)_\bbc \ \geq \ f(x_\e) + \left ({c / \sqrt n} \right )^n
 }
$$

The second inequality requires proof.  Abbreviate $x=x_\e$.  It suffices to show that each
point on the hypersurface $\Gamma$ defined by 
$\det_\bbc A_\bbc = (f^{1\over n}(x) + (c/\sqrt n))^n$ is distance $\leq c$ from the hypersurface
$\G'$ defined by $\det_\bbc B_\bbc = f(x)$.  Note that $B_0 \equiv f^{1\over n}(x) I \in \G'$
and the unit normal to $\G'$ at $B_0$ is $N={(1/ \sqrt n)}I$.  Set $A_0 = B_0 + cN
= ( f^{1\over n}(x) + (c / \sqrt n) ) I$.  Then $\det_\bbc A_0 =  ( f^{1\over n}(x) +  (c/\sqrt n) )^n$
so that $A_0\in\G$.  Note that $\dist(A_0, \G') =c$ and other points $A\in \G$ have smaller distance
to $\G'$.  Therefore, $\dist(A,\G') \geq c$ implies $\det_\bbc A \geq ( f^{1\over n}(x) +  (c/\sqrt n) )^n
\geq f(x) +  (c/\sqrt n )^n$.

Suppose now that case (i) of (2)$''$ holds, i.e., $-B_\e'\notin \cp_\bbc$.  
Then by (6)$'$(i) we have $\| {A_\e}'' + B_\e'\| = \dist\left({A_\e}'', -  B_\e'\ \right) \geq c$.
By (5),
$$
\| {A_\e}'' + B_\e'\|  \ =\ \L \|I\| {|x_\e-y_\e|^2 \over \e}  +  \left  \|(L_{x_e} - L_{y_\e})\left( {x_\e-y_\e  \over \e}\right)\right\|.
\eqno{(\HH.21)}
$$

This converges to zero as $\e\to 0$ by (3), so that this case (i) cannot occur for $\e>0$ small.

Finally, suppose Case (ii) of (2)$''$ holds.  By (6)$'$(ii) first and then (2)$''$(ii) we have:
$$
c^{1\over n} \ \leq\ \det_{\bbc}\left( {A_\e}'' \right)_\bbc - f(x_\e)
\ \leq\ 
\det_{\bbc}\left( {A_\e}'' \right)_\bbc - \det_{\bbc}\left( -  {B_\e}' \right)_\bbc
+f(y_\e) - f(x_\e),
\eqno{(\HH.22)}
$$
where  $( {A_\e}'' )_\bbc \geq0$ and $(-  {B_\e}' )_\bbc \geq0$.
Again by (5) and (3), as $\e\to0$, the RHS of (\HH.22) approaches zero, so this case cannot occur.\qed

      \vskip .3in
      
 \centerline{\bf Functions Which are Sub-the-Harmonics.}
 \medskip

      With regard to the subequation $F(J)$, since we have a notion of $F(J)$-harmonic functions
      (Definition \HH.2), one can also consider function that are ``sub''  these  harmonics.
      
      \Def{\HH.10}  A function $u\in\USC(X)$ is said to be {\bf sub-the-$F(J)$-harmonics}
      if for each compact set $K\ss X$ and each $F(J)$-harmonic function $h$ on a neighborhood
      of $K$,
      $$
      u\ \leq\ h \quad {\rm on}\ \ \partial K 
      \qquad\Rightarrow\qquad
      u\ \leq\ h \quad {\rm on}\ \  K. 
      $$
      
      Comparison for $F(J)$ on $K$ implies that each $F(J)$-subharmonic
      is also  sub-the-$F(J)$-harmonics.  Since $F(J)$-subharmonicity is a local property,
      the converse is an elementary consequence of local existence
      which follows as in Remark 9.7 and Theorem 9.2 in [HL$_4$]. This proves the following.
      
      \Theorem{\HH.11} {\sl
      Given a function $u\in\USC(X)$,
      \medskip
     \centerline
     {
     $u$ is $F(J)$-subharmonic 
         $\qquad \iff \qquad$
   $u$ is sub-the-$F(J)$-harmonics.
   }
}
     \medskip
      
      In the language of  [HL$_4$] this says that 
      $
      F(J)^{\rm visc} =      F(J)^{\rm classical}.
      $

      \vskip .3in
      
 \centerline{\bf $F$-Maximal Equals $F$-Harmonic.}
 \medskip
There is a difference  in the meaning of ``solution''   for our subequation $F=F(J,f)$,
depending on whether one is grounded in pluripotential theory or viscosity theory.

\Def {\HH.12} A function $u\in F(X)$ is  {\bf $F$-maximal} if  for each compact set $K\ss X$
\medskip
\centerline{$v\leq u$ on $\partial K \ \ \Rightarrow \ \ v\leq u$ on $K$}
\smallskip\noindent
for all $v\in \USC(K)$ which are $F$-subharmonic on $\Int K$.      
\medskip

The maximality property  for $u$ is equivalent to comparison holding for the function 
 $w\equiv -u$ (that is, $v+w$ satisfies the zero maximum principle for all functions $v$ which 
 are $F$-subharmonic).
 
 Now comparison holds for $F(J)$ and Perron functions $U$ are $F(J)$-harmonic on small balls.
 These facts can be used to show that following.
 
 \Prop{\HH.13} {\sl
 Given a function $u\in F(X)$,
\medskip
\centerline
{
$u$ is $F(J)$-harmonic $\IFF$ $u$ is $F(J)$-maximal.
 $$      
 }
 }

%%%%%%%%%%%%%%%%%%%%%%%%%%%%%%%%%%%%%%%%%%%%%%%
%%%%%%%%%%%%%%%%%%%%%%%%%%%%%%%%%%%%%%%%%%%%%%%
%%%%%%%%%%%%%%%%%%%%%%%%%%%%%%%%%%%%%%%%%%%%%%%
%%%%%%%%%%%%%%%%%%%%%%%%%%%%%%%%%%%%%%%%%%%%%%%
%%%%%%%%%%%%%%%%%%%%%%%%%%%%%%%%%%%%%%%%%%%%%%%

\vfill\eject
%\vskip .3in

\noindent{\headfont  \JJ.  
Distributionally Plurisubharmonic Functions.}\medskip

So far we have discussed three notions of plurisubharmonicity 
on an almost complex manifold: the viscosity notion, the standard notion
using restriction to holomorphic curves, and the notion of being sub-the-$F(J)$-harmonics.
In all three cases we start with the same object -- an upper semi-continuous function $u$.
Theorem \II.2 states that these first two notions are equal,
while Theorem \HH.7 states that the first and third notions are the same.  
There is a yet another definition of plurisubharmonicity
which starts with a distribution $u\in \cd'(X)$. (Let $v\geq 0$ stipulate that $v$ is a non-negative measure.)

\Def{\JJ.1}  A distribution $u\in \cd'(X)$ on an almost complex manifold 
$(X,J)$ is  {\bf distributionally $J$-plurisubharmonic on $X$} if
$$
\ch(V,V) (u) \ \geq\ 0 \fa V\in \G_{\rm cpt}(X, T_{1,0})
\eqno{(\JJ.1)}
$$
or equivalently
$$
H(v,v) (u) \ \geq\ 0  \fa v\in \G_{\rm cpt}(X, TX).
\eqno{(\JJ.2)}
$$

This distributional notion can not be ``the same'', 
but it is equivalent in a sense we now make precise.
In what follows we implicitly assume that an $\FJ$-plurisubharmonic function  
 $u\in\USC(X)$ is not  $\equiv -\infty$ on any component of $X$.

\Theorem {\JJ.2} {\sl \medskip
(a)\ \ Suppose $u$ is $\FJ$-plurisubharmonic.  Then $u\in \lloc(X) \ss \cd'(X)$, and 
$u$ is distributionally $J$-plurisubharmonic.

\medskip
(b) \ \  Suppose $u\in\cd'(X)$ is distributionally $J$-plurisubharmonic.  Then $u\in\lloc(X)$, and there exists a 
unique upper semi-continuous representative $\wt u$ of the $\lloc$-class $u$ which is  
$\FJ$-plurisubharmonic.  Moreover, $$\wt u (x) \ =\ {\rm ess} \limsup_{y\to x} u(y).$$ 
}

\medskip
\noindent
{\bf Remark.}  In light of Theorem \II.2, statement (a) above is equivalent to a result of Nefton Pali
[P, Thm. 3.9], and statement (b) provides a proof of his Conjecture 1 [P, p. 333].
Pali proves his conjecture under certain assumptions on the distrubution $u$ [P, Thm. 4.1].

\medskip
Both  notions of plurisubharmonicity in this theorem  can be reformulated using a family of ``Laplacians''.
To begin we consider the standard complex structure on $\bbc^n$.  Recall
(Example 4.4) that the standard subequation $\cp^\bbc\ss \Sym_\bbr(\bbc^n)$ is the set of real quadratic
forms (equivalently symmetric matrices) with non-negative $J_0$-hermitian symmetric part.  
The starting point is the following characterization of $\cp^\bbc$.
Assume $A\in  \Sym_\bbr(\bbc^n)$. Then
$$
A\ \in\ \cp^\bbc\qquad\iff\qquad \bra AB \ \geq\ 0 \quad\forall\, B \in H\Sym(\bbc^n)\ \ {\rm with}\ B>0.
\eqno{(\JJ.3)}
$$
The proof is left to the reader.

Each positive definite $B\in \Sym_\bbr(\bbr^N)$ defines a linear second-order operator
$$
\D_Bu\ =\ \bra{D^2u} B
\eqno{(\JJ.4)}
$$
which we call the $B$-Laplacian.

For $C^2$-fundtion $u$, the equivalence (\JJ.3) can be restated as a characterization of plurisubharmonic
functions.
$$
u\ \ {\rm is\ psh}   \qquad\iff\qquad \D_B u  \ \geq\ 0 \quad\forall\, B \in H\Sym(\bbc^n)\ \ {\rm with}\ B>0.
\eqno{(\JJ.3)'}
$$
Let
$$
H_B\ \equiv\  \{A\in  \Sym_\bbr(\bbc^n):  \bra AB \ \geq\ 0 \}
\eqno{(\JJ.5)}
$$
be the subequation associated to the differential operator $\D_B$.  Then (\JJ.3) can be restated as
$$
\cp^\bbc \ =\ \bigcap \left\{  H_B : 
 B \in H\Sym(\bbc^n)\  {\rm with}\ B>0 \right\}.
\eqno{(\JJ.3)''}
$$

Adopting the standard viscosity definition (using $C^2$-test functions) for $H_B$-subharmonic 
upper semi-continuous functions $u$, it is immediate from (\JJ.3)$''$  that for a $C^2$-function $\vf$,
which is a test function for $u$, we have
$$
 D^2_x \vf \in \cp^\bbc\qquad\iff\qquad 
  D^2_x \vf \in H_B \ \ \ \forall\,  B \in H\Sym(\bbc^n)\  {\rm with}\ B>0.
\eqno{(\JJ.6)}
$$
This proves the following.

\Prop{\JJ.3} {\sl
An upper semi-continuous function $u$ defined on an open subset $X$
of $\bbc^n$  is $\cp^\bbc$-plurisubharmonic if and only if it is $\D_B$-subharmonic for every
$B$-Laplacian, where $B \in H\Sym(\bbc^n)$ with  $B>0$.
}
\medskip

This proposition can be extended to $F(J)$-plurisubharmonic functions on any almost complex
manifold $(X,J)$, because of Proposition 4.5 (jet-equivalence).
Suppose that $g:X\to {\rm GL}^+_\bbr(\bbr^{2n})$ defines an almost complex structure
$J=g J_0 g^{-1}$ as in  Proposition 4.5, and let
 $E:X\to \Hom_\bbr(\bbc^n, \Sym_\bbr(\bbc^n))$
be defined as in (4.8).

\Def{\JJ.4}  Given  $B \in H\Sym(\bbc^n)$ with  $B>0$, define the $B$-Laplacian by
$$
L_B u \ =\ \bra {gBg^t}{D^2u}  + \bra  {E^t(gBg^t)}{Du}.
\eqno{(\JJ.7)}
$$

\Theorem{\JJ.5}  {\sl
An upper semi-continuous function $u$ on $X$ is $F(J)$-plurisubharmonic if and only if it is 
viscosity $L_B$-subharmonic for all  $B \in H\Sym(\bbc^n)$ with  $B>0$.
}

\pf  Apply (4.15) and Definition \JJ.4.\qed
\medskip

The parallel to Theorem \JJ.5 for distributional subharmonicity
is also true.  The proof is left to the reader.

\Def{\JJ.6}  A distribution $u\in \cd'(X)$ is {\bf distributionally $L_B$-subharmonic}
if $L_Bu$  is a non-negative measure on $X$.
\medskip

\Theorem{\JJ.7}  {\sl A distribution $u\in \cd'(X)$ on an almost complex manifold 
$(X,J)$ is   distributionally $J$-plurisubharmonic  if and only if 
locally, with $L_B$ defined by (\JJ.7), $u$ is distributionally $L_B$-subharmonic
for each  $B \in H\Sym(\bbc^n)$ with  $B>0$.}

\medskip

Combining Theorems \JJ.5 and \JJ.7 reduces Theorem \JJ.2 to the linear analogue for $L_B$.
This theorem is treated in the Appendix -- completing the proof of Theorem \JJ.2.\qed

 \Cor{\JJ.8}  {\sl  The concept of distributional $J$-plurisubharmonicity on an almost complex manifold $(X,J)$ is 
equivalent to the notion of standard $J$-plurisubharmonicity.}

\pf
Apply Theorem \JJ.2 and the restriction Theorem 6.1.\qed

%%%%%%%%%%%%%%%%%%%%%%%%%%%%%%%%%%%%%%%%%%%%%%%
%%%%%%%%%%%%%%%%%%%%%%%%%%%%%%%%%%%%%%%%%%%%%%%
%%%%%%%%%%%%%%%%%%%%%%%%%%%%%%%%%%%%%%%%%%%%%%%
%%%%%%%%%%%%%%%%%%%%%%%%%%%%%%%%%%%%%%%%%%%%%%%
%%%%%%%%%%%%%%%%%%%%%%%%%%%%%%%%%%%%%%%%%%%%%%%

%\vfill\eject
\vskip .3in

\noindent{\headfont  \GG. The Non-Equivalence of Hermitian and Standard Plurisubharmonic Functions.}\medskip

     Whenever an almost complex manifold $(X,J)$ is given a hermitian metric,  i.e., a riemannian
     metric $\bra \cdot\cdot$ such that $J_x$ is orthogonal for all $x$,  there is an induced notion
     of {\sl hermitian plurisubharmonicity} defined via the riemannian hessian (cf. [HL$_2$]).
     If the associated K\"ahler form $\o(v,w)=\bra {Jv} w$ is closed, then the hermitian
     plurisubharmonic functions agree with the intrinsic ones studied  in  this paper.
     However, in general they are not the same.  Proofs of these two assertions form the content
          of this section.
     
     To begin we recall that on a riemannian manifold $(X,\bra \cdot\cdot)$
     each smooth function $\vf$ has   a {\sl riemannian hessian} $\Hess \, \vf$ which is
     a section of $\Sym(T^*X)$  defined on (real) vector fields $v$ and $w$ by 
     $$
    (\Hess\, \vf)(v,w) \ =\ vw\vf -(\nabla_v w)\vf.
     $$
     where $\nabla$ denotes the Levi-Civita connection. If we are given $J$,  orthogonal
  with respect to $\bra \cdot\cdot$, then  
  $$
  (\Hess^\bbc \vf)(v,w) \ \equiv\  (\Hess\, \vf)(v,w) + (\Hess\, \vf)(Jv,Jw) 
 \eqno{(\GG.1)}
  $$
 is a well defined hermitian symmetric form on $(TX,J)$, i.e., a section of $\HB_\bbr(TX)$.
  A function  $\vf \in C^2(X)$ is then defined to be {\bf hermitian plurisubharmonic}
   if $\Hess^\bbc_x \vf \geq 0$  for all $x\in X$.
  This notion carries over to arbitrary upper semi-continuous functions,
  and these have been used to study complex Monge-Ampere equations
  in this setting (see [HL$_2$]).

     The natural question is:  How are these hermitian psh-functions related
     to the intrinsic ones?  In one important case they coincide.

     \Theorem{\GG.1} {\sl
     Suppose that $(X,J,\bra \cdot\cdot)$ is an almost complex hermitian manifold whose
     associated K\"ahler form is closed.  Then the space of hermitian plurisubharmonic
     functions on $X$ coincides with the space of standard plurisubharmonic functions.
     }

     \Remark{\GG.2}  Manifolds of this type are important in symplectic topology.
 Suppose that $(X,\o)$ is a given symplectic manifold and $J$ is an almost
 complex structure on $X$ such that
$$
\o(v,Jv) \ >\ 0\qquad \fa {\rm tangent\ vectors\ } v\neq 0.
$$
(Such $J$ always exist.) Then we can define an associated hermitian  metric by
$
\bra v w\ =\ \half\{\o(v,Jw) +\o(w,Jv)\}.
$
The triple $(X,J,\o)$ is sometimes called a {\sl Gromov manifold}.

\medskip

By a standard calibration argument, Theorem \GG.1 is a direct consequence of the following
more general result.

     \Theorem{\GG.3} {\sl
     Suppose that $(X,J,\bra \cdot\cdot)$ is an almost hermitian manifold in which
     all holomorphic curves are minimal surfaces (mean curvature 0).  
     Then the space of hermitian plurisubharmonic
     functions on $X$ coincides with the space of standard plurisubharmonic functions.
     }

\pf
By Definition \EE.8 we have 
$$
H(\vf)(v,v) \ =\  \{ v v   +(Jv)(Jv)  + J([v, Jv])\}\vf, 
\eqno{(\GG.2)}
  $$
which combined with Definition \GG.1 above gives the identity
$$
H(\vf)(v,v) \ =\  (\Hess^\bbc \vf)(v,v) + \{ \n_v v   +\n_{Jv}(Jv)  + J([v, Jv])\}\vf.
\eqno{(\GG.3)}
  $$
 Consider now the restriction of $H(\vf)$ to a holomorphic curve $\Sigma$, and suppose the vector field 
  $v$  is tangent to $\Sigma$ along $\Sigma$.  Then $Jv$ and $J[Jv,v]$ are
also tangent  to $\Sigma$ along $\Sigma$.  As we know  the restriction of $H(\vf)$ to 
$\Sigma$ agrees with the intrinsic hessian on $\Sigma$ and is given by the formula (\GG.2).

We now consider the restriction of the RHS of (\GG.3).  
The tangential part $\n^T_v v \equiv (\n_v v)^T$ is the intrinsic Levi-Civita
connection on $\Sigma$.  This connection is torsion-free and preserves $J$,
i.e., $\n^T J=0$.  (It preserves $J$ since it preserves the metric and $J$ is just rotation
by $\pi/2$.) Thus 
$$
\eqalign
{
 \n_v^T v   +\n_{Jv}^T(Jv)  + J([v, Jv]) \ &=\  \n_v^T v   +\n_{Jv}^T(Jv)  + J(\n_v^T(Jv) - \n_{Jv}^T v) \cr
 &=\  \n_v^T v   +\n_{Jv}^T(Jv)  + (\n_v^T(J^2v) - \n_{Jv}^TJ v) \cr
 &=\ 0.
 }
$$
Hence, equation (\GG.3) becomes
$$
H(\vf)(v,v) \ =\  (\Hess^\bbc \vf )(v,v) + \{ \n_v^N v   +\n_{Jv}^N(Jv) \}\vf
\eqno{(\GG.4)}
$$
 where $ \n_v^N v \equiv  \n_v v- \n_v^T v$.   This gives the following.

\Lemma {\GG.4} {\sl If $v$ is a vector field tangent to a holomorphic curve $\Sigma$,
then along $\Sigma$ we have
$$
H(\vf)(v,v) \ =\  (\Hess^\bbc \vf )(v,v) +  |v|^2 H_\Sigma \cdot \vf
\eqno{(\GG.5)}
$$
where $H_\Sigma$ is the mean curvature vector field of $\Sigma$ in the given hermitian metric.
}
\pf
By definition, the second fundamental form of $\Sigma$ is 
$B_{v,v} = \nabla^N_v v$.  Also by definition $H_\Sigma = \tr B =
B_{e,e} + B_{Je,Je}$ for a (any) unit tangent vector $e$ at the point. Formula (\GG.5)
now follows from (\GG.4).\qed
\medskip

Suppose now that all holomorphic curves are minimal.  Then since every tangent vector
$v$ on $X$ is tangent to a holomorphic curve by [NW], Lemma \GG.4 implies that
  the subequations
defined by $H(\vf)\geq0$ and $\Hess^\bbc(\vf)\geq0$ are the same.
This completes the proof of Theorem \GG.3. \qed

\vskip .3in

Now in general  things are not so nice.  Here are some  examples which show that the 
notions of hermitian and standard plurisubharmonicity are essentially unrelated.

\Ex{\GG.5}  Consider
$$
\bbc^2 \ = \ \bbr^4
$$
with coordinates $X=(z,w) = (x,y,u,v)$ and with the hermitian metric
$$
ds^2\ =\ {|dX|^2\over (1+|X|^2)^2}
\eqno{(\GG.6)}
$$
This is just the spherical metric in stereographic projection.
Consider the holomorphic curves
$$
\Sigma_x\ =\ \{(x,0)\}\times \bbc \ \ss\  \bbc\times\bbc
\eqno{(\GG.7)}
$$
Think of these on the 4-sphere.  They lie in the geodesic 3-sphere corresponding
to the 3-plane  $\{x$-axis $\}\times \bbc$. They form a family of round 2-spheres 
of constant mean curvature which are $\perp$ to the geodesic corresponding
to the x-axis, and which fill out the 3-sphere.  The mean curvature vector of $\Sigma_x$ is
$$
H^\Sigma_x\ =\ \phi(x){\partial\over \partial x}
$$ 
where
$$
x\phi(x) \ >\ 0 \ \ {\rm for}\ \  x\neq 0.
$$

Now consider the plurisubharmonic function
$$
u(z,w)\ =\ |f(z)|^2
$$
where $f$ is holomorphic. 
Since $u$ is constant on the curve $\Sigma_x$ we have $\Delta^\Sigma u=0$ and
by formula (\GG.5) we have
$$
\tr_{T_p\Sigma}\left\{\Hess\, u\right\} \ =\  - H^\Sigma\cdot u.
$$
Now choose $f(z)=z$ so that $u=|z|^2=x^2+y^2$.
Then
$$
\tr_{T_p\Sigma}\left\{\Hess\, u\right\} \ =\  - 2x\phi(x) \ <\ 0.
$$
We conclude that:
\medskip
\centerline{\sl Not every standard  plurisubharmonic function on $\bbc^2$}

\centerline{\sl
is hermitian plurisubharmonic for this hermitian metric.}
\medskip

Conversely, we now show that 

\medskip
\centerline{\sl Not every hermitian plurisubharmonic function on $\bbc^2$, with this metric,}

\centerline{\sl
is  plurisubharmonic in the standard sense.}
\medskip

Our example is local. To facilitate computation we make the following change of
coordinates.  Consider $S^4\ss\bbr^5$ (standardly embedded) and let 
$\Phi : S^4\to \bbr^4$ denote stereographic projection, where $\bbr^4=\bbc^2$ is our space above.
Fix a point $(a,0,0,0)\in\bbr^4$ ($a>0$).  Then there is a unique rotation $R$ of $S^4$ about the 
$(y,u,v)$-plane (i.e., a rotation of the $(x, \{{\rm vertical}\})$ 2-plane) which carries $\Phi^{-1}(a,0,0,0)$
to the south pole $(-1,0,0,0,0,)$.  Let $\Psi = \Phi\circ R\circ \Phi^{-1}$ be the change of coordinates.
Note that $\Psi(a,0,0,0)=(0,0,0,0)$.

Since $R$ is an isometry, the metric (\GG.6) is unchanged by this transformation.

Of course the complex structure on $\bbc^2$ has been conjugated.  In particular the holomorphic
curve $\Sigma_a$ in (\GG.7) has been transformed to the round 2-sphere passing through the origin:
$$
\Sigma \ \equiv\ \{(x,0,u,v) : (x-r)^2+u^2+v^2 = r^2\}
\eqno{(\GG.8)}
$$
where $r = r(a)>0$.

Consider now the function
$$
u(x,y,u,v) \ =\ \half\left( x^2+y^2+u^2+v^2\right) - Cx
$$
for $C>0$.  At the origin $0\in\bbr^4$ the hermitian hessian (= the riemannian 4-sphere hessian)
agrees with the standard coordinate hessian (since all the Christoffel symbols $\Gamma_{ij}^k$ vanish at 0).     
Thus
$$
\Hess_0 u = {\rm Id},
$$     
and we conclude that $u$ is hermitian plurisubharmonic in a neighborhood of 0.

On the other hand, the Laplace-Beltrami operator on $\Sigma$ satisfies
$$
\D_\Sigma u \ =\ \tr_{T\Sigma}( \Hess\, u) + H_\Sigma  \cdot u
$$
     where $H_\Sigma$ is the mean curvature vector of the 2-sphere $\Sigma$.
     Since $H_\Sigma = ( 2/r){\partial\over\partial x}$, we conclude that 
     $$
     (\D_\Sigma u )_0\ =\ 2 - {2\over r} C\ <\ 0
\eqno{(\GG.9)}
$$
if $C>r$.  The conformal structure on $\Sigma$ is the one induced from  $S^4$.  Hence (\GG.9) implies
that $u\bigr|_\Sigma$ is strictly superharmonic on a neighborhood of 0, and the assertion above is proved.

\medskip
\noindent
{\bf To Summarize:} {\sl  On the hermitian complex manifold $(\bbc^2, ds^2)$ not every standard
plurisubharmonic function is hermitian plurisubharmonic, and conversely,  not every hermitian
plurisubharmonic function is  plurisubharmonic in the standard sense.
}

\medskip
\noindent
{\bf Final Remark.}  The riemannian hessian is a bundle map 
$$
\Hess : J^2(X) \ \arr\ \Sym(T^*X)
\eqno{(\GG.10)}
$$
and, as such, splits the fundamental exact sequence
$
0 \to \Sym(T^*X) \to J^2(X) \to T^*X\to0.
$
The importance of this for nonlinear equations on riemannian manifolds
is the following.  Each  constant coefficient pure second-order subequation $\bbf_0\ss \Symn$,
which is O$(n)$-invariant, naturally induces a sub-fibre-bundle $F_0\ss \Sym(T^*X)$, 
and therefore determines  $F\ss J^2(X)$ by 
$$
F\ =\ \Hess^{-1}(F_0).
$$

The situation is analogous on an almost complex manifold $(X,J)$.  The real form of the 
complex hessian is a bundle map
$$
H : J^2(X) \ \arr\  {\rm HSym}^2(T^*X) \ \ss\ \Sym(T^*X).
\eqno{(\GG.11)}
$$
Any   $\bbf_0\ss {\rm HSym}^2(\bbc^n)$ which is GL$_\bbc(n)$-invariant and satisfies
$\bbf_0+\cp^\bbc \ss \bbf_0$ naturally induces a sub- sub-fibre-bundle $F_0\ss  {\rm HSym}^2(T^*X)$, 
and therefore determines  $F\ss J^2(X)$ by 
$$
F\ =\ H^{-1}(F_0).
$$

In both cases the resulting subequation is locally jet-equivalent to a constant coefficient equation.

There is ``overlap'' here corresponding to the cases Hess$^\bbc(\vf)\geq0$ and $H(\vf)\geq0$
discussed above. As shown in Example \GG.5, for general hermitian metrics these two
subequations are not the same.

 \medskip
 
 We note that the map (\GG.11) can be used to understand the polar cone 
 $F(J)^0 \ss J_2(X)$ (the lower 2-jet bundle).  It is the image of $\cp^\bbc \ss  {\rm HSym}^2(TX)$
 under the dual map $\Hess^*:  {\rm HSym}^2(TX) \to J_2(X)$.

%\vfill\eject
\vskip .3in

\centerline{\bf  References.}
\vskip .3in

\noindent
\item{[BT]}   E. Bedford and B. A. Taylor,  {The Dirichlet problem for a complex Monge-Amp\`ere equation}, 
Inventiones Math.{\bf 37} (1976), no.1, 1-44.

\smallskip

\item {[CKNS]} L. Caffarelli, J. J. Kohn,  L. Nirenberg, and J. Spruck, {\sl  The Dirichlet  problem for non-linear second order elliptic equations II: Complex Monge-Amp\`ere and uniformly elliptic equations},    Comm. on Pure and Applied Math. {\bf 38} (1985), 209-252.

\smallskip

\item {[HL$_{1}$]}  F. R. Harvey and H. B. Lawson, Jr., {\sl  Dirichlet duality and the non-linear Dirichlet problem},    Comm. on Pure and Applied Math. {\bf 62} (2009), 396-443.  ArXiv:0710.3991.

\smallskip

\item {[HL$_{2}$]}  ------,  {\sl   Dirichlet duality and the non-linear Dirichlet problem on Riemannian manifolds},  J. Diff. Geom.  {\bf 88} No. 3 (2011), 395-482.  ArXiv:0907.1981.

\smallskip

\item {[HL$_{3}$]}  ------, {\sl  The restriction theorem for fully nonlinear subequations},    
 {\sl Ann. Inst.  Fourier} (to appear).  ArXiv:1101.4850.

\smallskip

\item  {[HL$_{4}$]} \ \----------, {\sl The equivalence of  viscosity and distributional
subsolutions for convex subequations -- the strong Bellman principle},  Bol. Soc. Bras.  Mat.
(to appear).   ArXiv:1301.4914.

\smallskip

\item {[HL$_{5}$]}  \ \----------,   {\sl  Existence, uniqueness and removable singularities for nonlinear
partial differential equations in geometry}, 
Surveys in Geometry (to appear).  
\smallskip

\item {[NW]}  A. Nijenhuis and W. Woolf, {\sl  Some integration problems in almost -complex and complex manifolds}, 
Ann. of Math.,  {\bf 77} (1963), 424-489.

\smallskip

\item {[P]}  N. Pali, {\sl Fonctions plurisousharmoniques et courants positifs de type (1,1)
sur une vari\'et\'e presque complexe}, 
Manuscripta Math.  {\bf 118} (2005), no. 3, 311-337.

\smallskip

\item {[Pl]}  S. Pli\'s, {\sl The Monge-Amp\`ere equation on almost complex manifolds}, 
ArXiv:1106.3356, June, 2011  (v2, July, 2012).

%\vfill\eject
\vskip .5in

\centerline{\bf Appendix A.  The Equivalence of the Various Notions of Subharmonic}
\centerline{\bf for Reduced Linear Elliptic Equations }

\def\SH#1{{\rm SH}^{#1}}
\def\SHV{\SH {\rm visc}}
\def\SHW{\SH {\rm dist}}
\def\SHI{\SH {\rm class}}

\def \esssup#1{{\rm ess} \!\! \sup_{#1}}
\def \essup#1{{\rm ess}  \sup_{#1}}

\bigskip

  Consider the reduced linear operator %$L:C^\infty(X)\to C^\infty(X)$
%given by
$$
Lu(x) \ =\ a(x)\cdot D^2u(x) + b(x)\cdot Du(x)
$$
where $a$ and $b$ are $C^\infty$ on  an open set $X\ss \rn$, and $a>0$ is positive definite at each point.
A $C^2$-function $u$ on $X$ is said to be $L$-{\sl subharmonic} if $Lu\geq 0$ and $L$-{\sl harmonic} if
$Lu=0$. These  notions can be generalized using two completely different kinds of ``test functions''.

\Def {A.1}  A function $u\in \USC(X)$  is $L$-{\sl subharmonic in the viscosity sense} 
if for every $x\in X$ and every test function $\vf$ for $u$ at $x$, $L\vf(x)\geq0$  (cf. \S \CC).
Let $\SHV(X)$ denote the set of these.

\Def {A.2} A distribution $u\in \cd'(X)$ is  $L$-{\sl subharmonic in the distributional sense} if
$(Lu)(\vf) \equiv u(L^t\vf)\geq0$ for every non-negative test function $\vf\in C_{\rm cpt}^\infty(X)$, or equivalently, if
$Lu\geq 0$, i.e., $Lu$ is a positive measure on $X$.
Let $\SHW(X)$ denote the set of these.

\medskip

We say that $u$ is   {\sl viscosity $L$-harmonic} if both $u$ and $-u$ are viscosity $L$-subharmonic.
We say that $u$ is   {\sl distributionally $L$- harmonic} if  $Lu=0$ as a distribution.
In both cases there is a well developed theory of $L$-harmonics and $L$-subharmonics.

For example, the $L$-harmonics, both distributional and viscosity, are smooth.  This provides the proof 
that the two notions of $L$-harmonic are identical.  It  is not as straightforward to make statements relating
$\SHV(X)$ to $\SHW(X)$ since they are composed of different objects.
The bridge is partially provided by the following  third definition of $L$-subharmonicity.

\Def {A.3} A function  $u\in \USC(X)$  is {\sl classically $L$-subharmonic} 
if for every compact set $K\ss X$ and every $L$-harmonic function $\vf$ defined 
on a neighborhood of $K$, we have
$$
u\ \leq\ \vf \quad {\rm on} \ \partial K\qquad \Rightarrow
\qquad  u\ \leq\ \vf \quad {\rm on} \ K.
\eqno{(1)}
$$
Let $\SHI(X)$ denote the set of these.
\medskip

We always assume that $u$ is not identically $-\infty$ on any connected component of $X$.

In both the viscosity case and the distributional case a great number of results 
have been established including the following.

\Theorem{A.4} 
$$
\SHV(X)\ =\ \SHI(X)
$$

\Theorem{A.5} 
$$
\SHW(X)\ \cong\ \SHI(X)
$$

Note that Theorem A.4 can be stated as an equality since elements of both 
$\SHV(X)$ and $\SHI(X)$ are {\sl a priori} in $\USC(X)$.  By contrast, Theorem A.5
is not a precise statement until the isomorphism/equivalence is explicitly described.
%First we dispense with Theorem 4.

Theorem A.5 requires   careful attention, especially since in applications (such as the one in this paper)
the isomorphism sending $u\in \SHW(X)$ to $\wt u\in \SHI(X)$ is required to  produce the same upper
semi-continuous function $\wt u \in \USC(X)$ independent of the operator $L$.

 Separating out the two directions we have:  

\Theorem {A.6} 
{\sl
If $u\in \SHI(X)$, then $u\in \lloc(X) \ss \cd'(X)$ and $Lu\geq0$,  that is, $u\in \SHW(X)$.
} 
\medskip

\Theorem {A.7}
{\sl
Suppose $u\in \SHW(X)$.  Then $u\in \lloc(X)$.
Moreover, there exists 
 an upper semi-continuous representative  $ \wt u$ of the $\lloc$-class $u$ with $\wt u \in \SHI(X)$.
Furthermore,
$$
\wt u(x) \ \equiv \  \overrightarrow{ {\rm ess} \lim_{y\to x} }  u(y) \ \equiv \ \lim_{r\searrow0} {\rm  ess} \!\! \sup_{B_r(x)} u
$$
is the unique such representative.
}
\medskip

The precise statements, Theorem A.6 and Theorem A.7, give meaning to Theorem A.5.

Finally we outline some of the proofs.
\medskip
\noindent
{\bf Outline for Theorem A.4.}  
Assume $u\in \SHV(X)$ and $h$ is harmonic on a neighborhood of a compact set
$K\ss X$   with  $u\leq h$ on $\partial K$.
Since $\vf$ is a test function for $u$ at $x_0$ if an only if $\vf-h$ is  a test function for $u-h$ at $x_0$,
we have $u-h \in \SHV(X)$.  Since the maximum principle applies to $u-h$, we have $u\leq h$ on $K$.

Now suppose $u\notin \SHV(X)$.  Then there exists $x_0\in X$ and a test function
$\vf$ for $u$ at $x_0$ with $(L\vf)(x_0)<0$.  We can assume (cf. [HL$_2$, Prop. A.1]) that $\vf$ is a quadratic  and
$$
\eqalign
{
u-\vf\ &\leq\  -\l|x-x_0|^2 \quad {\rm for} \ |x-x_0|\leq \rho\ \ {\rm and}\cr
&= \ 0\qquad\qquad\qquad {\rm at} \ x_0
}
$$
for  some $\l,\rho>0$.  Set $\psi \equiv -\vf +\e$ where $\e= \l\rho^2$.
  Then $\psi$ is (strictly) subharmonic on a neighborhood of $x_0$.
Let $h$ denote the solution to the Dirichlet Problem on  $B\equiv B_{\rho}(x_0)$
with boundary values $\psi$. Since $h$ is the Perron function for $\psi\bigr|_{\partial B}$
and $\psi$ is $L$-subharmonic on $B$, we have $\psi\leq h$ on $\overline B$.
Hence, $-h(x_0) \leq -\psi(x_0)=\vf(x_0) -\e < u(x_0)$.  However,
on $\partial B$ we have $u \leq \vf -\l\rho^2 = -\psi  
=  -h$.  Hence, $u \notin \SHI(X)$.\qed

\medskip

\noindent
{\bf Outline for  Theorem A.6.}  
This theorem is part of classical potential theory, and a proof can be found in 
[HH], which also treats the hypo-elliptic case.
For fully elliptic operators $L$ we outline the part of the proof showing that $u\in\lloc(X)$.

Consider  $u\in \SHI(X)$. Fix a ball  $B\ss X$,
and let $P(x,y)$ be the Poisson kernel for the operator $L$ on $B$ (cf. [G]).
Then we claim that for $x\in\Int B$,
$$
u(x) \ \leq\ \int_{\partial B} P(x,y) u(y) d\s(y)
\eqno{(A.1)}
$$
where $\sigma$ is standard spherical measure.
To see this we first note that for $\vf \in C(\partial B)$, the
 unique solution to the Dirichlet problem for an $L$-harmonic function on $B$ with
boundary values $\vf$ is given by  $h(x) =  \int_{\partial B} P(x,y) \vf(y) d\s(y)$.
Since $u\in \SHI(X)$ we conclude that 
$$
u(x) \ \leq\ \int_{\partial B} P(x,y) \vf(y) d\s(y)
$$
for all $\vf   \in C(\partial B)$  with $u\bigr|_{\partial B}\leq \vf$.
The inequality (A.1) now follows since $u\bigr|_{\partial B}$ is u.s.c.,  and 
 $u\bigr|_{\partial B} = \inf\{\vf \in C(\partial B): u\leq \vf\}$.

Note that the integral (A.1)  is well defined (possibly $=-\infty$) since $u$ is bounded above.

%One now adapts  the arguments of H\"ormander [Hor, ???] to prove that $u\in\lloc(X)$.
Consider a family of concentric balls $B_r(x_0)$ in $X$ for $r_0\leq r\leq r_0+\kappa$
and suppose $x\in B_{r_0}$.  Then for any probability measure $\nu$ on the interval $[r_0, r_0+\kappa]$
we have 
$$
u(x) \ \leq\ \int_{[r_0, r_0+\kappa]} \int_{\partial B_r}  P_{r}(x,y) u(y) d\s(y) \, d\nu(r)
\eqno{(A.2)}
$$
where $P_r$ denotes the Poisson kernel for the ball $B_r$.
Let $E\ss X$ be the set of points $x$ such that $u$ is $L^1$ in a 
neighborhood of $x$.  Obviously $E$ is open.  Using (A.2) we conclude that
if $x\notin E$ then $u\equiv -\infty$ in a neighborhood of $x$ (cf.  [Ho$_1$,  Thm. 1.6.9]).  
Hence both $E$ and its complement
are open. Since we assume that $u$ is not $\equiv-\infty$ on any connected component of $X$,
we conclude that $u\in\lloc(X)$.

It remains to show that $Lu\geq 0$. This is exactly Theorem 1 on page 136 of [HH].\qed

\medskip

\noindent
{\bf Outline for Theorem A.7.}  In a neighborhood of any point $x_0\in X$ the distribution $u \in\SHW(X)$ 
is the sum of an $L$-harmonic function and a Green's potential
$$
v(x) \ =\ \int G(x,y) \mu(y)
\eqno {(A.3)}
$$
where $\mu\geq 0$ is a non-negative measure with compact support.
Here $G(x,y)$ is the Green's kernel for a ball $B$ about $x_0$.  It suffices to prove
Theorem A.7 for Green's potentials $v$ given by (A.3).  The fact that $v\in L^1(B)$ 
is a standard consequence of the fact  that $G \in L^1(B \times B)$ 
with singular support  on the diagonal. Since $G(x,y)\leq 0$, (A.3) defines a point-wise function
$v(x)$ near $x_0$ with values in $[-\infty, 0]$.  By replacing $G(x,y)$ with the continuous kernel
$G_n(x,y)$, defined to be the maximum of $G(x,y)$ and $-n$, the integrals 
$v_n(x) = \int G_n(x,y)\mu(y)$ provide a decreasing sequence of continuous functions converging
to $v$.  Hence, $v$ is upper semi-continuous.  The maximum principle states that $v\in \SHI(X)$.

Finally we prove that if $u\in \lloc(X)$ has a representative $v\in \SHI(X)$, then $v=\wt u$.
Since 
$$
\esssup {B_r(x)} u \ =\ \esssup {B_r(x)} v \ \leq \ \sup_{B_r(x)} v,
\eqno{(A.4)}
$$
and $v$ is upper semi-continuous, it follows that $\wt u(x) \leq v(x)$. 

Applying (A.2) to $v$ and using the fact that $\int P_r(x,y) d\s(y) =1$  (since 1 is $L$-harmonic)  yields
 $$
 \eqalign
 {
v(x_0) \ &\leq\   {1\over \kappa} \int_{[0, \kappa]} \int_{\partial B_r} P_{r}(x_0,y) v(y) d\sigma(y) \, dr \cr
&\leq \ \left ( \essup{B_\kappa} v\right)   {1\over \kappa}  \int_{[0, \kappa]}  \int_{\partial B_r}  P_{r}(x_0,y)   d\sigma(y) \, dr \cr
& = \ \essup{B_\kappa} v \ =\ \essup {B_\kappa} u  \cr
}
$$
proving that $v(x_0) \leq \wt u(x_0)$.
\qed

\Remark{A.8}
The construction of $\wt u$ above is quite general and enjoys several nice properties, which we include here.
To any function $u\in \lloc(X)$ we can associate its {\sl essential upper semi-continuous regularization}:
$$
\wt u(x) \ =\ \overrightarrow{ {\rm ess} \lim_{y\to x} } u(y) \ \equiv \ \lim_{r\searrow0} {\rm  ess} \!\! \sup_{B_r(x)} u
$$
This clearly depends only on the $\lloc$-class of $u$.

\Lemma {A.9} {\sl For any $u\in \lloc(X)$, the function $\wt u$ is upper semi-continuous. 
Furthermore, for any $v\in \USC(X)$ representing the $\lloc$-class  $u$, we have   $\wt u\leq  v$,
and if $x\in X$ is a Lebesgue point for $u$ with value $u(x)$,
then $u(x)\leq \wt u(x)$.}

\pf
To show that $\wt u$ is upper semi-continuous, i.e., $\limsup_{y\to x} \wt u(y) \leq \wt u(x)$, it 
suffices to show that  
$$
\sup_{B_r(x)} \wt u  \leq \esssup {B_r(x)} u
$$ 
and then let $r\searrow 0$.  
However, if $B_\rho(y) \ss B_r(x)$, then 
$$
\wt u (y)\ =\ \lim_{\rho\to0} \esssup {B_\rho(y)} u\leq \esssup  {B_r(x)} u.
$$

Letting $r\searrow0$ in (A.4) proves that $\wt u (x)\leq v(x)$.

For the last assertion  of the lemma suppose that $x$ is a Lebesgue point for $u$ with value $u(x)$, i.e., by definition
$$
\lim_{r\to0} {1\over \left| B_r(x)\right|} \int |u(y)-u(x)| \, dy \ =\ 0,
$$
and hence the value $u(x)$ must be the limit of the means
$$
u(x) \ =\ \lim_{r\to0} {1\over \left| B_r(x)\right|} \int u(y) \, dy \ \leq\ \lim_{r\to0} \esssup{B_r(x)} u \ =\ \wt u(x).
\qquad \mathqed
$$

\vfill\eject

\Remark {A.10}  Theorem A.4 holds for any linear operator $L$ with $a\geq0$ provided that
\medskip

$\bullet$ \ \ $L$-harmonic functions are smooth,

\medskip

$\bullet$ \ \ The maximum principle holds for $u\in \SHV(X)$

\medskip

$\bullet$ \ \ The Perron function gives the unique solution for the Dirichlet problem for $L$.

\medskip
\noindent
{\bf Some Historical/Background Remarks.}
The equivalence of $\SHV(X)$ and $\SHW(X)$ for linear elliptic operators has been addressed by Ishii [I],
who proves the result for continuous functions but leaves open the case where 
$u\in\SHV(X)$ is a general upper semi-continuous function and the case where
$u\in \SHW(X)$ is a general distribution.
The proof that ``classical implies distributional'' appears in [HH] where the result is proved for 
even more general linear hypoelliptic operators $L$. Other arguments that ``viscosity implies distributional'' 
are known to Hitoshi Ishii and to Andrzej Swiech.
A treatment of mean value characterizations of $L$-subharmonic functions (again for subelliptic $L$)
can be found in [BL]. A general introduction to viscosity theory appears in [CIL].
A good discussion of the Greens kernel  appears in ([G]), and  the explicit  construction of the Hadamard parametrix
is found in  [Ho$_2$, 17.4]. 

 \medskip
 
 The arguments outlined here are presented in more detail and for general 
 convex subequations on manifolds in [HL].

%\vfill\eject
\vskip .3in

\centerline{\bf  References for Appendix A.}
\vskip .3in

\item {[BL]}   A. Bonfiglioli and E. Lanconelli  {\sl  Subharmonic functions in sub-riemannian settings},   to appear in the 
Journal of the European Mathematical Society.

\medskip

\item {[CIL]} M. G. Crandall, H. Ishii and P.-L. Lions, 
{\sl  User's guide to viscosity solutions of second order
 partial differential equations},    Bull. A. M. S.  {\bf 27} (1992), 1-67.

\smallskip

\item {[G]} P. Garabedian, { Partial Differential Equations},    J. Wiley and Sons,  New York, 1964.

\smallskip

\item {[HH]}  M. Herv\'e and R.M. Herv\'e., {\sl  Les fonctions surharmoniques  dans l'axiomatique de
M. Brelot associ\'ees \`a un op\'erateur elliptique d\'eg\'en\'er\'e},    Annals de l'institut Fourier,  {\bf 22}, no. 2 (1972), 131-145.

\item  {[HL]} \ \----------, {\sl The equivalence of  viscosity and distributional
subsolutions for convex subequations -- the strong Bellman principle},  Bol. Soc. Bras. de Mat. (to appear). ArXiv:1301.4914.

\smallskip

\smallskip

\item {[Ho$_1$]} L. H\"ormander, {An Introduction to Complex Analysis in Several Variables},    D. Von Nostrand Co., 
Princeton, 1966.

\smallskip

\item {[Ho$_2$]} L. H\"ormander, {The Analysis of Linear Partial Differential Operators, III},    Springer-Verlag, 
Berlin Heidelberg, 1985.

\smallskip

\item {[I]} H. Ishii, {\sl  On the equivalence of two notions of weak solutions, viscosity solutions and distribution solutions},    Funkcialaj Ekvacioj, {\bf 38} (1995), 101-120.

         %%%%%%%%%%%%% 
    \end